\documentclass[12pt]{article}
\usepackage{a4wide}

\usepackage{a4,epsfig}
\usepackage{amsfonts}
\usepackage{pgfplots}
\usepackage{pgfplotstable}
\usepackage{booktabs}
\usepackage{float}
\usepackage{hyperref}
\usepackage{amsmath}

\newtheorem{theorem}{Theorem}[section]

\newtheorem{corollary}[theorem]{Corollary}
\newtheorem{lemma}[theorem]{Lemma}

\newtheorem{remark}{Remark}[section]

\newcommand{\myproof} [1]
   { \noindent {\bf Proof.} #1 \hfill\rule{0.5em}{1.2ex} \par\medskip}

\newcommand{\nnorm}[1]{%
  \mid\!\mid\!\mid \! #1 \! 
  \mid\!\mid\!\mid}

\numberwithin{equation}{section} 

\begin{document}

\setcounter{page}{1}

\title{Space-time tensor-product finite element methods \\[1mm]
  for parabolic problems}
\author{Richard L\"oscher$^1$, Michael Reichelt$^2$, Olaf~Steinbach$^1$}
\date{$^1$Institut f\"ur Angewandte Mathematik, TU Graz, \\
  Steyrergasse 30, 8010 Graz, Austria \\[2mm]
  $^2$Lehrstuhl f\"ur Numerische Mathematik, 
  Helmut--Schmidt--Universit\"at, \\[1mm]
  Holstenhofweg 85, 22043 Hamburg, Germany}

\maketitle

\begin{abstract}
 We study space-time Galerkin--Petrov formulations for parabolic evolution
 problems and their relation to classical implicit time-stepping schemes.
 Although such schemes are stable in the usual time-stepping sense, their
 interpretation as space-time operator equations may lead to conditional
 stability, with constants depending on the relation between temporal and
 spatial mesh sizes. We revisit this phenomenon for the continuous Galerkin
 method of Aziz and Monk, which yields the Crank--Nicolson scheme in the
 lowest-order case, and provide a detailed space-time error analysis for
 solutions of both high and low regularity. In particular, the space-time
 framework allows us to analyze the deteriorated behaviour of classical
 time-stepping methods for nonsmooth initial data. By applying integration
 by parts in time, we derive an adjoint space-time formulation that
 incorporates the initial condition in a natural variational way. In the
 lowest-order case, this formulation leads to a Rannacher-type smoothing
 of the initial data. The theoretical results are complemented by
 numerical experiments.

\end{abstract}

\section{Introduction}
In this work, we are interested in the numerical solution of parabolic
evolution equations using space-time methods. 
As a model problem, we consider the initial Dirichlet boundary value
problem for the heat equation in a bounded Lipschitz domain
$\Omega \subset {\mathbb{R}}^n$, $n=1,2,3$, and a finite time horizon $T>0$,
\begin{equation}\label{Heat equation}
  \begin{array}{rclcl}
    \partial_t u(x,t) - \Delta_x u(x,t)
    & = & f(x,t)
    && \mbox{for} \; (x,t) \in Q := \Omega \times (0,T), \\[1mm]
    u(x,t)
    & = & 0
    && \mbox{for} \; (x,t) \in \Sigma := \partial \Omega \times (0,T), \\[1mm]
    u(x,0) & = & u_0(x) && \mbox{for} \; x \in \Omega .
  \end{array}
\end{equation}
A standard approach for the numerical solution of \eqref{Heat equation} is to
first discretize the problem in space, e.g. by a finite element or finite
difference method, and then solve the resulting system of ordinary differential
equations by a time stepping scheme. The analysis of such schemes is well
established, see, e.g., \cite{HairerNorsettWanner:1993}. However, a classical
time stepping analysis relies in most cases on point values of the solution and
the right-hand side in time, called temporal collocation. Consequently, the
resulting error estimates are usually formulated pointwise in time and
require more regularity of the data than would be needed for estimates in
natural Sobolev or Bochner spaces.

A way to overcome this limitation is to consider Galerkin--Petrov variational
formulations in time. These methods can be divided into two main classes: In
continuous Galerkin (cG) methods, the ansatz space consists of continuous
functions in time, whereas the test space is discontinuous. Such schemes were
studied, for instance, in \cite{AzizMonk:1989,Hulme:1971,Wright:1969}.
In contrast, in discontinuous Galerkin (DG) methods, both
ansatz and test spaces are discontinuous in time,
see \cite{ErikssonJohnsonThomee:1985,Lesaint:1987,LesaintRaviart:1974}.
Depending on the quadrature formula used for the temporal integrals, these
Galerkin--Petrov schemes are equivalent to well known time stepping methods.
At the same time, their variational structure provides a framework for
deriving error estimates under weaker regularity assumptions.

Instead of applying a Galerkin--Petrov formulation only to the purely temporal
problem, one may also consider space-time variational formulations. These
methods make use of the fact that the space-time operator equation is
well-posed in suitable space-time Banach spaces. In such approaches, space
and time are treated simultaneously, leading to a monolithic system rather
than to a sequential time stepping procedure, see, for instance,
\cite{DieningStevensonStorn:2025,KoetheLoescherSteinbach:2026,Steinbach:2015,
  SteinbachYang:2019,StevensonWesterdiep:2021,SchwabStevenson:2009,Ta:2025}.
This point of view is particularly attractive for
adaptive methods \cite{Akrivis:2011,Feischl:2022}, time-parallel algorithms
\cite{Gander:2015}, and the construction of preconditioners for the resulting
monolithic systems \cite{DieningStevensonStorn:2025,NeumuellerSmears:2019},
as well as optimal control problems subject to time dependent partial differential
equations \cite{DanielsHinzeVierling:2015}.
Space-time discretizations may either be built on general unstructured
space-time meshes or on tensor product spaces. In both cases, stability of
the discrete space-time formulation is a central issue.

A detailed analysis for cG schemes on tensor product meshes was carried out
by Andreev \cite{Andreev:2012}. A central observation of this analysis is that,
although the corresponding time stepping scheme may be unconditionally stable
in the classical sense, the associated space-time operator equation can suffer
from conditional stability. More precisely, the stability constant may depend
on the relation between the temporal and spatial mesh sizes, leading to a
CFL-type restriction that enforces parabolic scaling. Similar conclusions were
obtained in \cite{DieningStorn:2022,DieningStevensonStorn:2025}. Thus, the
space-time perspective does not merely reproduce classical time stepping
schemes, but also provides a variational framework for analyzing their
deteriorated behaviour for nonsmooth data. In particular, the loss of accuracy
observed for low-regularity initial conditions in classical time stepping
methods \cite{Rannacher:1985,Osterby2003} can be understood at the level of
the associated space-time operator equation, see also
\cite{AndreevSchwitzer:2014}.

In this paper, we focus first on the cG scheme introduced by Aziz and Monk
\cite{AzizMonk:1989}, which yields the Crank--Nicolson method in the
lowest order case. We revisit the conditional stability result of
\cite{Andreev:2012} in a space-time framework and analyze the resulting
discretization by means of an inf-sup argument. In particular, we give a
detailed space-time error analysis and derive convergence rates for solutions
of both high and low regularity. Furthermore, by applying integration by parts
in time, we introduce an adjoint space-time formulation. This formulation
incorporates the initial condition in a natural way and, in the lowest order
case, leads to a scheme that can be interpreted as a Rannacher-type
smoothing of the initial data \cite{Rannacher:1985}.

The paper is structured as follows. In Section~\ref{sec:primal continuous}, we
introduce the space-time variational formulation for the model problem and
recall existence and uniqueness results based on the inf-sup theory. In
Section~\ref{sec:error estimates space time}, we collect space-time
approximation estimates for functions of different regularity, which are then
used in Section~\ref{sec:primal discrete} to derive error estimates for
the primal formulation. In Section~\ref{sec:adjoint discrete}, we apply
integration by parts in time and obtain the corresponding adjoint formulation.
The analysis is complemented by numerical experiments in
Section~\ref{sec:numerics}. Finally, conclusions are given in
Section~\ref{sec:Conclusions}.
  
\section{The primal space-time variational formulation}
\label{sec:primal continuous}
The primal space-time variational formulation of the initial boundary
value problem \eqref{Heat equation} is to find $u \in X$ with
$u(\cdot,0)=u_0 \in L^2(\Omega)$, such that
\begin{equation}\label{Heat equation VF}
  \int_0^T \int_\Omega \Big[
  \partial_t u(x,t) \, v(x,t) + \nabla_x u(x,t) \cdot
  \nabla_x v(x,t) \Big] dx \, dt \, = \,
  \int_0^T \int_\Omega f(x,t) \, v(x,t) \, dx \, dt
\end{equation}
is satisfied for all $v \in Y := L^2(0,T;H^1_0(\Omega))$, where
$X := \{ u \in Y : \partial_t u \in Y^* \}$. Related norms are given as
\[
  \| v \|_Y := \| \nabla_x v \|_{L^2(Q)}, \quad
  \| u \|_X := \sqrt{\| u \|_Y^2 + \| \partial_t u \|^2_{Y^*}}
  = \sqrt{\| u \|^2_Y + \| w_u \|^2_Y},
\]
where $w_u \in Y$ is the unique solution of the variational formulation
\begin{equation*}
  \int_0^T \int_\Omega \nabla_x w_u(x,t) \cdot
  \nabla_x v(x,t) \, dx \, dt \, = \, \int_0^T \int_\Omega
  \partial_t u(x,t) \, v(x,t) \, dx \, dt
\end{equation*}
for all $v \in Y$. For $u \in X$ and $v \in Y$ we introduce the bilinear
form
\begin{equation}\label{eq:bilinearform b}
  b(u,v) \, := \,
  \int_0^T \int_\Omega \Big[
  \partial_t u(x,t) \, v(x,t) + \nabla_x u(x,t) \cdot
  \nabla_x v(x,t) \Big] dx \, dt,
\end{equation}
satisfying
\begin{equation}\label{Bound b}
  |b(u,v)| \leq \sqrt{2} \, \| u \|_X \| v \|_Y \quad
  \mbox{for all} \; (u,v) \in X \times Y .
\end{equation}
To ensure the well-posedness of the variational formulation
\eqref{Heat equation VF}, we assume $f \in Y^*$, and $u_0 \in L^2(\Omega)$.
The latter ensures the
existence of an extension $\widetilde{u}_0 \in X$ satisfying
$\widetilde{u}_0(\cdot,0)=u_0\in L^2(\Omega)$. It remains to determine
$\widetilde{u} \in X_0 := \{ u \in X : u(\cdot,0)=0 \in L^2(\Omega) \}$
as unique solution of the variational problem
\begin{equation}\label{Heat equation VF hom}
  b(\widetilde{u}+\widetilde{u}_0,v) \, = \, \langle f , v \rangle_Q
  \quad \mbox{for all} \; v \in Y .
\end{equation}
Unique solvability of the variational formulation
\eqref{Heat equation VF hom} follows from the surjectivity of the
bounded bilinear form $b(\cdot,\cdot)$, and from
the inf-sup stability condition
\begin{equation}\label{Heat inf sup}
  \| u \|_X \leq \sup\limits_{0 \neq v \in Y}
  \frac{b(u,v)}{\| v \|_Y} \quad \mbox{for all} \; u \in X_0 .
\end{equation}
The proof of \eqref{Heat inf sup} follows when considering the
particular test function $\overline{v} := u + w_u$ for $u \in X_0 \subset Y$,
see \cite[Theorem 2.1]{Steinbach:2015}.
For conforming space-time finite element spaces $X_{0,h} \subset X_0$ and
$Y_h \subset Y$, the Galerkin variational formulation of
\eqref{Heat equation VF hom} is to find $\widetilde{u}_h \in X_{0,h}$ such that
\begin{equation}\label{Heat equation FEM}
  b(\widetilde{u}_h+\widetilde{u}_0,v_h) \, = \, \langle f , v_h \rangle_Q
  \quad \mbox{for all} \; v_h \in Y_h .
\end{equation}
Unique solvability of \eqref{Heat equation FEM} and related space-time
finite element error estimates follow from the discrete inf-sup
stability condition
\begin{equation}\label{Heat inf sup discrete}
  \| u_h \|_{X,h} \leq \sup\limits_{0 \neq v_h \in Y_h}
  \frac{b(u_h,v_h)}{\| v_h \|_Y} \quad \mbox{for all} \; u_h \in X_{0,h},
\end{equation}
when using the discrete norm
\[
\| u \|_{X,h} := \sqrt{\| u \|_Y^2 + \| w_{u,h} \|^2_Y},
\]
where $w_{u,h} \in Y_h$ solves the variational formulation
\begin{equation*}\label{Heat equation Def wh}
  \langle \nabla_x w_{u,h} , \nabla_x v_h \rangle_{L^2(Q)} =
  \langle \partial_t u , v_h \rangle_Q \quad
  \mbox{for all} \; v_h \in Y_h .
\end{equation*}
Indeed, the discrete inf-sup stability condition \eqref{Heat inf sup discrete}
follows as in the continuous case when choosing
$\overline{v}_h := u_h + w_{u_h,h} \in Y_h$
for $u_h \in X_{0,h} \subset Y_h$. In \cite{Steinbach:2015},
we have considered space-time finite element spaces
$X_{0,h}=Y_h=\mbox{span} \{ \varphi_k \}_{k=1}^M  \subset X_0 \subset Y$
of piecewise linear continuous basis functions $\varphi_k$ which are
defined with respect to an admissible locally quasi-uniform decomposition
of the space-time domain $Q$ into simplicial shape regular space-time
finite elements $q_\ell$ of mesh size $h_\ell$. 
However, in this paper we will consider different tensor product
space-time finite element spaces $X_{0,h} \subset X_0$ and $Y_h \subset Y$ 
satisfying $\mbox{dim} \, X_{0,h} = \mbox{dim} \, Y_h$
but $X_{0,h} \not\subset Y_h$, in order
to reformulate the space-time Galerkin formulation
\eqref{Heat equation FEM} in the spirit of more common
time stepping schemes. For simplicity, in this paper we only
consider lowest order space-time finite element spaces, but this
approach can be extended to higher order polynomial spaces as well.

\section{Space-time tensor product finite element spaces}
\label{sec:error estimates space time}
For the spatial domain $\Omega \subset {\mathbb{R}}^n$ we consider an
admissible and globally quasi-uniform decomposition into shape regular
simplicial finite elements $\tau_\ell$ of mesh size $h_x$, and the
related finite element space
$V_{h_x} = \mbox{span} \{ \varphi_k \}_{k=1}^{M_x} \subset H^1_0(\Omega)$
of piecewise linear continuous basis functions $\varphi_k$.
In the case of a globally quasi-uniform finite element mesh
there holds the inverse inequality
\begin{equation}\label{inverse inequality}
  \| \nabla_x u_h \|_{L^2(\Omega)} \, \leq \, c \, h_x^{-1} \,
  \| u_h \|_{L^2(\Omega)} \quad \mbox{for all} \; u_h \in V_{h_x}.
\end{equation}
For $u \in L^2(\Omega)$ we introduce the $L^2$ projection
$Q_{h_x}^1u \in V_{h_x}$ as the unique solution of the variational formulation
\begin{equation*}\label{Def Qhx}
  \langle Q_{h_x}^1 u , v_h \rangle_{L^2(\Omega)} \, = \,
  \langle u , v_h \rangle_{L^2(\Omega)} \quad
  \mbox{for all} \; v_h \in V_{h_x},
\end{equation*}
for which we immediately conclude the trivial stability estimate
\begin{equation}\label{Stability L2}
  \| Q_{h_x}^1 u \|_{L^2(\Omega)} \, \leq \, \| u \|_{L^2(\Omega)} \quad
  \mbox{for all} \; u \in L^2(\Omega).
\end{equation}
When assuming $u \in H^1_0(\Omega) \cap H^s(\Omega)$ for some $s \in [1,2]$, 
there hold the following error estimates, e.g.,
\cite{BrennerScott:1994,Steinbach:2008}:
\begin{eqnarray}
  \| u - Q_{h_x}^1 u \|_{L^2(\Omega)} \label{Error Qx L2}
  & \leq & c \, h_x^s \, |u|_{H^s(\Omega)} , \\
  \| \nabla_x (u - Q_{h_x}^1 u) \|_{L^2(\Omega)} \label{Error Qx H1}
  & \leq & c \, h_x^{s-1} \, |u|_{H^s(\Omega)} . 
\end{eqnarray}
Using \eqref{Error Qx L2} and \eqref{Error Qx H1}, and a space interpolation
argument for some $\tau \in [0,1]$, we conclude the error estimate
\begin{equation}\label{Error Qx Htau}
  \| u - Q_{h_x}^1 u \|_{H^\tau(\Omega)} \, \leq \,
  c \, h_x^{s-\tau} \, |u|_{H^s(\Omega)} .
\end{equation}
Moreover, e.g., \cite{BramblePasciakSteinbach:2002}, 
$Q_{h_x}^1 : H^1_0(\Omega) \subset L^2(\Omega) \to
V_{h_x} \subset H^1_0(\Omega)$ is stable in $H^1_0(\Omega)$, i.e.,
\begin{equation}\label{H1 Stability}
  \| \nabla_x Q_{h_x}^1u \|_{L^2(\Omega)} \, \leq \,
  c_Q \, \| \nabla_x u \|_{L^2(\Omega)} \quad
  \mbox{for all} \; u \in H^1_0(\Omega).
\end{equation}
In addition to the $L^2$ projection $Q_{h_x}^1u$ we introduce,
for $ u \in H^1_0(\Omega)$, the $H^1_0$ projection $P_{h_x}^1u \in V_{h_x}$
as the unique solution of the variational formulation
\begin{equation}\label{Def Phx}
  \langle \nabla_x P_{h_x}^1 u , \nabla_x v_h \rangle_{L^2(\Omega)} \, = \,
  \langle \nabla_x u , \nabla_x v_h \rangle_{L^2(\Omega)} \quad
  \mbox{for all} \; v_h \in V_{h_x},
\end{equation}
where we conclude the stability estimate
\begin{equation}\label{Stability Ph}
  \| \nabla_x P_{h_x}^1 u \|_{L^2(\Omega)} \, \leq \,
  \| \nabla_x u \|_{L^2(\Omega)} \quad \mbox{for all} \; u \in H^1_0(\Omega),
\end{equation}
and the error estimates
\begin{eqnarray}
  \| \nabla_x (u - P_{h_x}^1 u) \|_{L^2(\Omega)} \label{Error Px H1}
  & \leq & c \, h_x^{s-1} \, |u|_{H^s(\Omega)} , \\ 
  \| u - P_{h_x}^1 u \|_{L^2(\Omega)} \label{Error Px L2}
  & \leq & c \, h_x^s \, |u|_{H^s(\Omega)} , 
\end{eqnarray}
when assuming $u \in H^1_0(\Omega) \cap H^s(\Omega)$ for some $s \in [1,2]$.

For the time interval $(0,T)$ we consider a uniform decomposition into
$N_t$ finite elements $(t_{i-1},t_i)$ of mesh size $h_t = T/N_t$, and
nodes $t_i=ih_t$ for $i=0,\ldots,N_t$. With respect to this temporal mesh
we introduce the space
$S_{h_t}^1(0,T) = \mbox{span} \{ \psi_i^1 \}_{i=0}^{N_t}$ of
piecewise linear and continuous basis functions $\psi_i^1$. For continuous
$ u \in H^\sigma(0,T)$,
$\sigma > \frac{1}{2}$, we define the nodal interpolation
\begin{equation}\label{Def Iht}
(I_{h_t}^1u)(t) = \sum\limits_{i=0}^{N_t} u(t_i) \psi_i^1(t) \, ,
\end{equation}
for which we have the error estimate
\begin{equation}\label{Error It}
  \| u - I_{h_t}^1 u \|_{L^2(0,T)} \leq c \, h_t^\sigma \,
  |u|_{H^\sigma(0,T)},
\end{equation}
when assuming $u \in H^\sigma(0,T)$ for some $\sigma \in (1/2,2]$. Note that
the constant degenerates as $\sigma \to 1/2$. Finally,
$W_{h_t}^0 = \mbox{span} \{ \psi_i^0 \}_{i=1}^{N_t} \subset L^2(0,T)$
is the space of piecewise constant basis functions $\psi_i^0$, for which the
temporal $L^2$ projection $Q_{h_t}^0u \in W_{h_t}^0$ is defined as the unique
solution of the variational formulation
\begin{equation}\label{Def Qht}
  \langle Q_{h_t}^0 u , v_h \rangle_{L^2(0,T)} =
  \langle u , v_h \rangle_{L^2(0,T)} \quad \mbox{for all} \; v_h \in W_{h_t}^0,
\end{equation}
satisfying
the error estimate
\begin{equation}\label{Error Qht}
  \| u - Q_{h_t}^0 u \|_{L^2(0,T)} \, \leq \,
  c \, h_t^\sigma \, |u|_{H^\sigma(0,T)},
\end{equation}
when assuming $u \in H^\sigma(0,T)$ for some $\sigma \in [0,1]$.
For $ t \in (t_{i-1},t_i)$ we finally compute
\[
  \partial_t I_{h_t}^1 u(t) = \frac{1}{h_t} \, (u(t_i)-u(t_{i-1})) =
  \frac{1}{h_t} \int_{t_{i-1}}^{t_i} \partial_t u(t) \, dt =
  Q_{h_t}^0\partial_t u(t) \, .
\]
Now we are in a position to state some approximation error estimates
in space-time tensor product finite element spaces:

\begin{lemma}\label{Lemma Approximation QhxIht}
  Let $u \in L^2(0,T;H^1_0(\Omega) \cap H^s(\Omega)) \cap
  H^\sigma(0,T;H^1_0(\Omega))$ for some $s \in [1,2]$ and $\sigma \in (1/2,2]$.
  Then there holds the error estimate
  \begin{equation}\label{Error QI}
    \| \nabla_x(u - Q_{h_x}^1 I_{h_t}^1 u) \|_{L^2(Q)} \, \leq \,
    c_1 \, h_x^{s-1} \, |u|_{L^2(0,T;H^s(\Omega))} +
    c_2 \, h_t^\sigma \, |u|_{H^\sigma(0,T;H^1_0(\Omega))} .
  \end{equation}  
  If $u \in L^2(0,T;H^1_0(\Omega) \cap H^s(\Omega)) \cap
  H^\sigma(0,T;L^2(\Omega))$ is satisfied for some $s \in [1,2]$ and
  $\sigma \in (1/2,2]$, then  there holds the error estimate
  \begin{equation}\label{Error QI 2}
    \| \nabla_x (u - Q_{h_x}^1 I_{h_t}^1u) \|_{L^2(Q)} \, \leq \,
    c_1 \, h_x^{s-1} \, |u|_{L^2(0,T;H^s(\Omega))} +
    c_2 \, h_x^{-1} \, h_t^\sigma \, |u|_{H^\sigma(0,T;L^2(\Omega))} .
  \end{equation}
\end{lemma}
  
\myproof{By the triangle inequality, we first have
  \begin{equation}\label{Step 1}
    \| \nabla_x (u - Q_{h_x}^1 I_{h_t}^1 u) \|_{L^2(Q)} \, \leq \,
    \| \nabla_x (u - Q_{h_x}^1 u) \|_{L^2(Q)} +
    \| \nabla_x Q_{h_x}^1 (u - I_{h_t}^1 u) \|_{L^2(Q)} .
  \end{equation}
  For the first part we can use the error estimate (\ref{Error Qx H1})
  for $s \in [1,2]$ to conclude
  \begin{eqnarray*}
    \| \nabla_x (u - Q_{h_x}^1 u) \|_{L^2(Q)}^2
    & = &
    \int_0^T \| \nabla_x (u(t) - Q_{h_x}^1 u(t)) \|_{L^2(\Omega)}^2 \, dt
    \\ & \leq &
    c \, h_x^{2(s-1)} \, \int_0^T |u(t)|^2_{H^s(\Omega)} \, dt,
  \end{eqnarray*}
  i.e.,
  \[
    \| \nabla_x (u - Q_{h_x}^1 u) \|_{L^2(Q)}
    \, \leq \, c \, h_x^{s-1} \, |u|_{L^2(0,T;H^s(\Omega))} .
  \]
  To bound the second term in \eqref{Step 1}, we first use the stability
  estimate \eqref{H1 Stability} and the interpolation error estimate
  \eqref{Error It} to obtain
  \begin{eqnarray*}
    && \| \nabla_x Q_{h_x}^1 (u - I_{h_t}^1 u) \|_{L^2(Q)}^2 \, = \, \int_0^T
       \| \nabla_x Q_{h_x}^1(u(t)-I_{h_t}^1u(t)) \|_{L^2(\Omega)}^2 \, dt \\
    && \hspace*{2mm} \leq \, c_Q^2 \, \int_0^T
       \| \nabla_x (u(t)-I_{h_t}^1u(t)) \|_{L^2(\Omega)}^2 \, dt \, = \,
       c_Q^2 \, \int_\Omega \| (I - I_{h_t}^1) \nabla_x u(x)
       \|^2_{L^2(0,T)} \, dx \\
    & & \hspace*{2mm} \leq \, c \, h_t^{2\sigma} \, \int_\Omega
        \| \nabla_x u \|^2_{H^\sigma(0,T)} \, dx \, = \,
        c \, h_t^{2\sigma} \, \| \nabla_x u \|^2_{H^\sigma(0,T;L^2(\Omega))},
  \end{eqnarray*}
  i.e.,
  \[
    \| \nabla_x Q_{h_x}^1 (u - I_{h_t}^1 u) \|_{L^2(Q)} \, \leq \,
    c \, h_t^\sigma \, \| u \|_{H^\sigma(0,T;H^1_0(\Omega))} ,
  \]
  and \eqref{Error QI} follows.
      
  However, when using the inverse inequality \eqref{inverse inequality}
  in space first, the stability estimate \eqref{Stability L2}, and the
  temporal interpolation error estimate \eqref{Error It}, we obtain
  \begin{eqnarray*}
    \| \nabla_x Q_{h_x}^1 (u - I_{h_t}^1 u) \|_{L^2(Q)}^2
    & = & \int_0^T \| \nabla_x Q_{h_x}^1 (u - I_{h_t}^1 u)
          \|_{L^2(\Omega)}^2 \, dt \\  
    & \leq & c \, h_x^{-2} \,
             \int_0^T \| Q_{h_x}^1 (u(t)-I^1_{h_t}u(t)) \|^2_{L^2(\Omega)} \, dt \\
    & \leq & c \, h_x^{-2} \,
             \int_0^T \| u(t)-I_{h_t}^1u(t) \|^2_{L^2(\Omega)} \, dt \\
    & = & c \, h_x^{-2} \, \int_\Omega
          \| u(x) - I^1_{h_t} u(x) \|^2_{L^2(0,T)} \, dx \\
    & \leq & c \, h_x^{-2} \, h_t^{2\sigma} \, \int_\Omega
             |u(x)|^2_{H^\sigma(0,T)} \, dx ,
  \end{eqnarray*}
  i.e.,
  \[
    \| \nabla_x Q_{h_x}^1 (u - I_{h_t}^1 u) \|_{L^2(Q)} \, \leq \,
    c \, h_x^{-1} \, h_t^\sigma \, \| u \|_{H^\sigma(0,T;L^2(\Omega))} ,
  \]
  and this gives \eqref{Error QI 2}.}

\noindent
When the nodal piecewise linear interpolation $I_{h_t}^1$ in time is
replaced by the temporal piecewise constant $L^2$ projection $Q_{h_t}^0$,
the statement of Lemma \ref{Lemma Approximation QhxIht} can be reformulated
as follows, where the proof can be adapted accordingly.

\begin{lemma}\label{Lemma Approximation QhxQht}
  Let $u \in L^2(0,T;H^1_0(\Omega) \cap H^s(\Omega)) \cap
  H^\sigma(0,T;H^1_0(\Omega))$ for some $s \in [1,2]$ and $\sigma \in [0,1]$.
  Then there holds the error estimate
  \begin{equation}\label{Error QQ}
    \| \nabla_x(u - Q_{h_x}^1 Q_{h_t}^0 u) \|_{L^2(Q)} \, \leq \,
    c_1 \, h_x^{s-1} \, |u|_{L^2(0,T;H^s(\Omega))} +
    c_2 \, h_t^\sigma \, |u|_{H^\sigma(0,T;H^1_0(\Omega))} .
  \end{equation}  
  If $u \in L^2(0,T;H^1_0(\Omega) \cap H^s(\Omega)) \cap
  H^\sigma(0,T;L^2(\Omega))$ is satisfied for some $s \in [1,2]$ and
  $\sigma \in [0,1]$, then  there holds the error estimate
  \begin{equation}\label{Error QQ 2}
    \| \nabla_x (u - Q_{h_x}^1 Q_{h_t}^0 u) \|_{L^2(Q)} \, \leq \,
    c_1 \, h_x^{s-1} \, |u|_{L^2(0,T;H^s(\Omega))} +
    c_2 \, h_x^{-1} \, h_t^\sigma \, |u|_{H^\sigma(0,T;L^2(\Omega))} .
  \end{equation}
\end{lemma}

\noindent
It remains to state an approximation result for $Q_{h_x}^1 I_{h_t}^1$ in
$Y^*$:

\begin{lemma}\label{Lemma Approximation Y*}
  For $u \in X\cap H^\varrho(0,T;H^{-1}(\Omega))$, let $\partial_t u \in L^2(0,T;H^{-\tau}(\Omega))$ be satisfied for some $\tau \in [0,1]$,
  and $\varrho \in [1,2]$. Then there holds the error estimate
  \begin{equation}\label{Error QI t}
    \| \partial_t (u - Q_{h_x}^1 I_{h_t}^1 u) \|_{Y^*} \, \leq \,
    c_1 \, h_x^{1-\tau} \, \| \partial_t u \|_{L^2(0,T;H^{-\tau}(\Omega))}
    \, + \, c_2 \, h_t^{\varrho-1} \, \| u \|_{H^\varrho(0,T;H^{-1}(\Omega))} .
  \end{equation}
\end{lemma}
    
\myproof{For any $v \in Y = L^2(0,T;H^1_0(\Omega))$, we first consider
  \begin{eqnarray*}
    \langle \partial_t (u-Q_{h_x}^1I_{h_t}^1 u) , v \rangle_Q
    & = & \langle \partial_t u - Q_{h_x}^1 \partial_t u , v \rangle_Q +
          \langle Q_{h_x}^1(\partial_t u - \partial_t I_{h_t}^1 u),v
          \rangle_Q \\
    & = & \langle \partial_t u , (I-Q_{h_x}^1)v \rangle_Q +
          \langle \partial_t u - \partial_t I_{h_t}^1 u , Q_{h_x}^1 v
          \rangle_Q \, .
  \end{eqnarray*}
  For the first part, we obtain, using duality and the error
  estimate \eqref{Error Qx Htau} for some $\tau \in (0,1]$,
  \begin{eqnarray*}
    \langle \partial_t u , (I-Q_{h_x}^1)v \rangle_Q
    & \leq & \| \partial_t u \|_{L^2(0,T;H^{-\tau}(\Omega))} \,
             \| (I - Q_{h_x}^1) v \|_{L^2(0,T;H^\tau_0(\Omega))} \\
    & \leq & c \, h_x^{1-\tau} \,
             \| \partial_t u \|_{L^2(0,T;H^{-\tau}(\Omega))} \,
             \| v \|_{L^2(0,T;H^1_0(\Omega))} \, . 
  \end{eqnarray*}
  For the remaining term we conclude, using
  $\partial_t I_{h_t}^1 u = Q_{h_t}^0\partial_t u$, duality,
  the stability estimate \eqref{H1 Stability}, and the error
  estimate \eqref{Error Qht},
  \begin{eqnarray*}
    \langle \partial_t u - \partial_t I_{h_t}^1 u , Q_{h_x}^1 v \rangle_Q
    & = & \langle (I-Q_{h_t}^0) \partial_t u , Q_{h_x}^1 v \rangle_Q \\
    & \leq & \| (I-Q_{h_t}^0) \partial_t u \|_{L^2(0,T;H^{-1}(\Omega))} \,
             \| Q_{h_x}^1 v \|_{L^2(0,T;H^1_0(\Omega))} \\
    & \leq & c \, h_t^{\varrho-1} \,
             \| u \|_{H^\varrho(0,T;H^{-1}(\Omega))} \,
             \| v \|_{L^2(0,T;H^1_0(\Omega))},
  \end{eqnarray*}
  and the assertion follows from the norm definition of
  $\| \partial_t (u-Q_{h_x}^1 I_{h_t}^1u) \|_{Y^*}$ by duality.}

\section{The primal space-time finite element method}
\label{sec:primal discrete}
For a space-time tensor product finite element formulation of
\eqref{Heat equation FEM} we consider the ansatz space
$X_{0,h} := W^1_{h_t} \otimes V_{h_x} \subset X_0$ with
$W_{h_t}^1 = S_h^1(0,T) \cap H^1_{0,}(0,T) 
= \mbox{span} \{ \psi_i^1 \}_{i=1}^{N_t} $, and the test space
$Y_h := W_{h_t}^0 \otimes V_{h_x} \subset Y$. By construction we have
$\mbox{dim} \, X_{0,h} = \mbox{dim} \, Y_h = N_t \, M_x$, but
$X_{0,h} \not \subset Y_h$. In this case, unique solvability of
\eqref{Heat equation FEM} is based on the following result:

\begin{lemma}
  For $X_{0,h} = W^1_{h_t} \otimes V_{h_x}$ and $Y_h = W_{h_t}^0 \otimes V_{h_x}$
  there holds the discrete inf-sup condition
  \begin{equation}\label{Heat equation inf sup L2}
    \| \partial_t u_h \|_{L^2(Q)} \, \leq \,
    \sup\limits_{0 \neq v_h \in Y_h} \frac{b(u_h,v_h)}{\| v_h \|_{L^2(Q)}}
    \quad \mbox{for all} \; u_h \in X_{0,h} .
  \end{equation}
\end{lemma}
\myproof{For $u_h \in X_{0,h} = W_{h_t}^1 \otimes V_{h_x}$, we have
  $\overline{v}_h := \partial_t u_h \in Y_h = W_{h_t}^0 \otimes V_{h_x}$,
  and hence we can write
\begin{eqnarray*}
  b(u_h,\overline{v}_h)
  & = & \langle \partial_t u_h , \partial_t u_h \rangle_{L^2(Q)} +
        \langle \nabla_x u_h , \nabla_x \partial_t u_h \rangle_{L^2(Q)} \\[2mm]
  & = & \| \partial_t u_h \|^2_{L^2(Q)} +
        \int_0^T \int_\Omega \nabla_x u_h(x,t) \cdot
        \nabla_x \partial_t u_h(x,t) \, dx \, dt \\[2mm]
  & = & \| \partial_t u_h \|_{L^2(Q)}^2 +
        \frac{1}{2} \int_0^T \frac{d}{dt} \int_\Omega
        |\nabla_x u_h(x,t)|^2 \, dx \, dt \\[2mm]
  & = & \| \partial_t u_h \|_{L^2(Q)}^2 +
        \frac{1}{2} \, \| \nabla_x u_h(T) \|^2_{L^2(\Omega)} \\
  & \geq &
        \| \partial_t u_h \|_{L^2(Q)}^2 \, = \,
        \| \partial_t u_h \|_{L^2(Q)}
        \| \overline{v}_h \|_{L^2(Q)},
\end{eqnarray*}
and the assertion follows.}

\noindent
While we can use \eqref{Heat equation inf sup L2} to ensure unique
solvability of the variational formulation \eqref{Heat equation FEM}, we
cannot use \eqref{Heat equation inf sup L2} to derive a priori space-time
finite element error estimates for the numerical solution
$\widetilde{u}_h \in X_{0,h}$. Hence, we will state a different discrete
inf-sup condition as follows: When using the temporal piecewise constant
$L^2$ projection $Q_{h_t}^0$ as defined in \eqref{Def Qht}, we introduce the
mesh dependent norm, similar to the one used in \cite{UrbanPatera:2012},
\begin{equation}
  \nnorm{u}_X \, := \, \Big[ c_Q^{-2} \, \| \partial_t u \|_{Y^*}^2 + 
  \| Q^0_{h_t} u \|_Y^2 \Big]^{1/2} ,
  \label{upper norm equivalence} 
\end{equation}
where $c_Q$ is the stability constant as used in \eqref{H1 Stability}.

\begin{lemma}
  For $X_{0,h} = W^1_{h_t} \otimes V_{h_x}$ and $Y_h = W_{h_t}^0 \otimes V_{h_x}$,
  there holds the discrete inf-sup condition
  \begin{equation}\label{Heat equation inf sup Error}
    \nnorm{u_h}_X
    \, \leq \,
    \sup\limits_{0 \neq v_h \in Y_h} \frac{b(u_h,v_h)}{\| v_h \|_Y}
    \quad \mbox{for all} \; u_h \in X_{0,h} .
  \end{equation}
\end{lemma}
\myproof{Let $u_h\in X_{0,h} = W_{h_t}^1 \otimes V_{h_x}$ be arbitrary but fixed.
  For $\partial_t u_h \in W_{h_t}^0 \otimes V_{h_x} = Y_h$ we determine
  $w_{u_h,h}\in Y_h$ as the unique solution of the variational formulation
\begin{equation}\label{eq:auxiliary w_h for u_h}
  \langle \nabla_x w_{u_h,h},\nabla_x v_h \rangle_{L^2(Q)} =
  \langle \partial_t u_h,v_h \rangle_{L^2(Q)} \quad
  \mbox{for all} \; v_h\in Y_h. 
\end{equation}
It holds that 
\begin{eqnarray*}
  \| \nabla_x w_{u_h,h} \|_{L^2(Q)}^2
  & = & \langle \nabla_x w_{u_h,h},\nabla_x w_{u_h,h} \rangle_{L^2(Q)} \\
  & = & \langle \partial_t u_h,w_{u_h,h} \rangle_{L^2(Q)} \, \leq \,
        \| \partial_t u_h \|_{Y^*} \| \nabla_x w_{u_h,h} \|_{L^2(Q)}, 
\end{eqnarray*}
from which $\| w_{u_h,h} \|_Y\leq \| \partial_t u_h \|_{Y^*}$ follows.
In order to prove the inverse estimate, we use the definitions and
stability estimates for the $L^2$ projections $Q_{h_x}^1$ and $Q_{h_t}^0$,
respectively, to conclude, recall $\partial_t u_h \in Y_h$,
\begin{eqnarray*}
  \| \partial_t u_h \|_{Y^*}
  & = & \sup_{0 \neq v\in Y}
        \frac{\langle \partial_t u_h,v\rangle_{L^2(Q)}}{\| v \|_Y}
        \, = \,
        \sup_{0 \neq v\in Y}
        \frac{\langle \partial_t u_h,Q_{h_t}^0Q_{h_x}^1v \rangle_{L^2(Q)}}
        {\| v \|_Y} \\
  & = & \sup_{0\neq v\in Y} \frac{\langle \nabla_x w_{u_h,h} ,
        \nabla_x Q_{h_t}^0Q_{h_x}^1v  \rangle_{L^2(Q)}}{\| v \|_Y} \\
  & \leq & \sup_{0\neq v\in Y} \frac{\| \nabla_x w_{u_h,h} \|_{L^2(Q)}
        \| \nabla_x Q_{h_t}^0Q_{h_x}^1v  \|_{L^2(Q)}}{\| v \|_Y}
        \, \leq \,
        c_Q \, \| w_{u_h,h} \|_Y .
\end{eqnarray*}
For $\overline{v}_h := Q_{h_t}^0 u_h + w_{u_h,h} \in Y_h$ we now compute,
using \eqref{eq:auxiliary w_h for u_h}, 
\begin{eqnarray*}
  \| \overline{v}_h \|_Y^2
  & = & \| Q_{h_t}^0 u_h + w_{u_h,h} \|_Y^2 \\[1mm]
  & = & \| Q_{h_t}^0 u_h \|_Y^2 + \| w_{u_h,h} \|_Y^2 +
        2 \, \langle w_{u_h,h} , Q_{h_t}^0 u_h \rangle_Y \\[1mm]
  & \geq & \| Q_{h_t}^0 u_h \|_Y^2 + c_Q^{-2} \, \| \partial_t u_h \|_{Y^*}^2
           + 2 \, \langle \nabla_x w_{u_h,h}, \nabla_x Q_{h_t}^0 u_h
           \rangle_{L^2(Q)} \\[1mm]
  & = & \nnorm{u_h}_X^2 + 2 \, \langle \partial_t u_h , Q_{h_t}^0u_h
        \rangle_{L^2(Q)} \\[1mm]
  & = & \nnorm{u_h}_X^2 + 2 \,
        \langle \partial_t u_h, u_h \rangle_{L^2(Q)} \, \geq \,
        \nnorm{u_h}_X^2,
\end{eqnarray*}
due to
\[
  \langle \partial_t u_h,u_h \rangle_{L^2(Q)} =
  \int_\Omega \int_0^T \frac{1}{2} \frac{d}{dt} |u_h(x,t)|^2 \, dt \, dx
  = \frac{1}{2} \int_\Omega |u_h(x,T)|^2 \, dx \geq 0 .
\]
Moreover, we have that 
\begin{eqnarray*}
  && b(u_h,\overline{v}_h)
  \, = \, \langle \partial_t u_h,Q_{h_t}^0 u_h + w_{u_h,h} \rangle_{L^2(Q)} +
        \langle \nabla_x u_h , \nabla_x (Q_{h_t}^0 u_h + w_{u_h,h})
        \rangle_{L^2(Q)} \\[1mm]
  && = \, \langle \nabla_x w_{u_h,h},
        \nabla_x(Q_{h_t}^0 u_h + w_{u_h,h}) \rangle_{L^2(Q)} +
        \langle \nabla_x Q_{h_t}^0 u_h , \nabla_x (Q_{h_t}^0 u_h + w_{u_h,h})
        \rangle_{L^2(Q)} \\[1mm]
  && = \, \| \nabla_x (Q_{h_t}^0u_h+w_{u_h,h}) \|^2_{L^2(Q)} \, = \,
        \| \overline{v}_h \|_Y^2 \, \geq \,
        \| \overline{v}_h \|_Y \nnorm{u_h}_X,
\end{eqnarray*}
which gives the desired result.}

\noindent
The stability of the $L^2$ projection in $H^1_0(\Omega)$ was a sufficient
ingredient of the proof. In fact, one can show that it is also necessary,
see \cite{TantardiniVeeser:2016}. We are now in the position to prove a
first best approximation result. 

\begin{lemma}
  Let $\widetilde{u}\in X_0$ and $\widetilde{u}_h\in X_{0,h}$ denote the
  unique solutions of \eqref{Heat equation VF hom}
  and \eqref{Heat equation FEM}, respectively. Then there holds
  Cea's lemma,
  \begin{equation}\label{eq:Ceas lemma primal}
    \nnorm{\widetilde{u}-\widetilde{u}_h}_X \, \leq \,
    \sqrt{2} \, c_Q \,
    \inf\limits_{z_h \in X_{0,h}} \| \widetilde{u} - z_h \|_X .
  \end{equation}
\end{lemma}

\myproof{Let $G_h:X_0\to X_{0,h}$ denote the Galerkin projection, defined as 
  \[
    b(G_h u,v_h) = b(u,v_h),\quad \forall v_h\in Y_h, \; u \in X_0. 
  \]
  With the discrete inf-sup condition \eqref{Heat equation inf sup Error}
  we then conclude
  \[
    \nnorm{G_h u}_X \, \leq \,
    \sup_{0\neq v_h\in Y_h} \frac{b(G_hu,v_h)}{\| v_h \|_Y} \, = \,
    \sup_{0\neq v_h\in Y_h} \frac{b(u,v_h)}{\| v_h \|_Y}
    \, \leq \, \sqrt{2} \, c_Q \, \nnorm{u}_X,
  \]
  where we have used
  \begin{eqnarray*}
    b(u,v_h)
    & = & \langle \partial_t u , v_h \rangle_Q +
          \langle \nabla_x u , \nabla_x v_h \rangle_{L^2(Q)} \\[1mm]
    & = & \langle \partial_t u , v_h \rangle_Q +
          \langle \nabla_x Q_{h_t}^0 u , \nabla_x v_h \rangle_{L^2(Q)} \\[1mm]
    & \leq & \| \partial_t u \|_{Y^*} \| v_h \|_Y +
             \| Q_{h_t}^0 u \|_Y \| v_h \|_Y \\
    & \leq & \sqrt{2} \, \Big[ \| \partial_t u \|^2_{Y^*} +
             \| Q_{h_t}^0 u \|_Y^2 \Big]^{1/2} \| v_h \|_Y \, \leq \,
             \sqrt{2} \, c_Q \, \nnorm{u}_X \| v_h \|_Y .
  \end{eqnarray*}}

\noindent
Moreover, using that $G_h z_h = z_h$ for any $z_h\in X_{0,h}$,
$\nnorm{u}_X \leq \|u\|_X$, and
$\nnorm{ (I-G_h) } = \nnorm{ G_h }$, see \cite{Carstensen:2025}, we have 
  \begin{eqnarray*}
    \nnorm{\widetilde{u}- \widetilde{u}_h}_X
    & = & \nnorm{(I-G_h)\widetilde{u}}_X \, = \,
          \nnorm{(I-G_h)(\widetilde{u}-z_h)}_X \\[1mm]
    & \leq & \nnorm{ I-G_h } \, \nnorm{ \widetilde{u}-z_h }_X \, \leq \,
             \sqrt{2} \, \| \widetilde{u}-z_h \|_X. 
  \end{eqnarray*}

\noindent
When combining \eqref{eq:Ceas lemma primal} with the approximation
properties of the space-time finite element space $X_{0,h}$ as given
in Lemma \ref{Lemma Approximation QhxIht} and
Lemma \ref{Lemma Approximation Y*}, we then
conclude the following result:

\begin{lemma}
  Let $\widetilde{u}\in X_0$ and $\widetilde{u}_h\in X_{0,h}$ denote the
  unique solutions of \eqref{Heat equation VF hom}
  and \eqref{Heat equation FEM}, respectively. Assume
  $\widetilde{u} \in H^{2,1}(Q) \cap H^2(0,T;H^{-1}(\Omega))$ such that
  $\nabla_x \partial_t \widetilde{u} \in L^2(Q)$ is satisfied, and let
  $h_t \simeq h_x$. Then there holds the error estimate
  \begin{eqnarray}\label{Error primal smooth}
    \nnorm{\widetilde{u} - \widetilde{u}_h}_X
    & \leq & c \, h_x \, \Big[
             |\widetilde{u}|_{L^2(0,T;H^2(\Omega))}^2
             +
             \| \widetilde{u} \|_{H^1(0,T;L^2(\Omega))}^2 \\
    && \nonumber \hspace*{15mm}
             +
             \| \nabla_x \partial_t \widetilde{u} \|_{L^2(Q)}^2
             +
             \| \widetilde{u} \|_{H^2(0,T;H^{-1}(\Omega))}^2 \Big]^{1/2} .
  \end{eqnarray}
\end{lemma}

\myproof{We first consider the error estimate \eqref{Error QI} for $s=2$ and
  $\sigma=1$,
  \[
    \| \widetilde{u} - Q_{h_x}^1 I_{h_t}^1 \widetilde{u} \|_Y
    \, \leq \, c_1 \, h_x \, |\widetilde{u}|_{L^2(0,T;H^2(\Omega))} +
    c_2 \, h_t \, | \widetilde{u} |_{H^1(0,T;H^1_0(\Omega))} \, ,
  \]
  and the error estimate \eqref{Error QI t} for $\tau=0$ and
  $\varrho = 2$, i.e.,
  \[
    \| \partial_t (\widetilde{u} - Q_{h_x}^1 I_{h_t}^1\widetilde{u}) \|_{Y^*}
    \, \leq \, c_1 \, h_x \, 
    \| \widetilde{u} \|_{H^1(0,T;L^2(\Omega))} +
    c_2 \, h_t \, \| \widetilde{u} \|_{H^2(0,T;H^{-1}(\Omega))} . 
  \]
  Together with \eqref{eq:Ceas lemma primal} for
  $z_h = Q_{h_x}^1 I_{h_t}^1 \widetilde{u}$ and 
  choosing $h_t \simeq h_x$ we obtain the assertion.}

\noindent
Before we can prove an error estimate for
$\| \widetilde{u} - \widetilde{u}_h \|_X$, we need to establish
an error estimate for $\| \widetilde{u} - \widetilde{u}_h \|_{L^2(Q)}$, as
\eqref{Error primal smooth} only involves the discrete norm
$\| Q_{h_t}^0(\widetilde{u} - \widetilde{u}_h) \|_Y$:

\begin{lemma}\label{Lemma L2 Error primal smooth}
  Let $\widetilde{u}\in X_0$ and $\widetilde{u}_h\in X_{0,h}$ denote the
  unique solutions of \eqref{Heat equation VF hom}
  and \eqref{Heat equation FEM}, respectively. Assume
  $\widetilde{u} \in 
  L^2(0,T;H^2(\Omega)) \cap H^2(0,T;L^2(\Omega))$ such that
  $\partial_t \widetilde{u} \in L^2(0,T;H^2(\Omega))$ and
  $\widetilde{u} \in H^2(0,T;H^2(\Omega))$ is satisfied.
  For $h_t \simeq h_x$ there holds the error estimate
  \begin{equation}\label{L2 Error primal smooth} 
    \| \widetilde{u} - \widetilde{u}_h \|_{L^2(Q)} \, \leq \, c \, h_x^2 \,
    \Big[
    \| \widetilde{u} \|_{H^1(0,T;H^2(\Omega))} \, + \,
    \| \widetilde{u} \|_{H^2(0,T;H^2(\Omega))} \Big] \, .
   \end{equation}
\end{lemma}

\myproof{When using the spatial $H^1_0$ projection $P_{h_x}^1 : H^1_0(\Omega)
  \to V_{h_x}$ as defined in \eqref{Def Phx}, 
  and the piecewise linear temporal
  interpolation operator $I_{h_t}^1 : H^1_{0,}(0,T) \to W_{h_t}^1$,
  see \eqref{Def Iht}, we can write, using the Galerkin orthogonality
  \begin{equation}\label{Proof L2 smooth Galerkin}
    \langle \partial_t \widetilde{u}_h , v_h \rangle_Q +
    \langle \nabla_x \widetilde{u}_h , \nabla_x v_h \rangle_{L^2(Q)}
    =
     \langle \partial_t \widetilde{u} , v_h \rangle_Q +
    \langle \nabla_x \widetilde{u} , \nabla_x v_h \rangle_{L^2(Q)}
  \end{equation}
  for the particular test function, 
  $v_h = \partial_t (\widetilde{u}_h - I_{h_t}^1 P_{h_x}^1 \widetilde{u})
  \in Y_h$, and using
  $\partial_t I_{h_t}^1 \widetilde{u} = Q_{h_t}^0 \partial_t \widetilde{u} $,
  \begin{eqnarray*}
    \| \partial_t (\widetilde{u}_h - I_{h_t}^1 P_{h_x}^1 \widetilde{u})
    \|^2_{L^2(Q)}
    & = & \langle \partial_t (\widetilde{u}_h - I_{h_t}^1 P_{h_x}^1
          \widetilde{u}),
          \partial_t (\widetilde{u}_h - I_{h_t}^1 P_{h_x}^1
          \widetilde{u}) \rangle_{L^2(Q)} \\[1mm]
    &&  \hspace*{-4cm}
       \leq \, \langle \partial_t (\widetilde{u}_h - I_{h_t}^1 P_{h_x}^1
             \widetilde{u}),
             \partial_t (\widetilde{u}_h - I_{h_t}^1 P_{h_x}^1 \widetilde{u})
             \rangle_{L^2(Q)} \\[1mm]
    && + \,
       \langle \nabla_x (\widetilde{u}_h - I_{h_t}^1 P_{h_x}^1 \widetilde{u}),
       \nabla_x \partial_t (\widetilde{u}_h - I_{h_t}^1 P_{h_x}^1
       \widetilde{u}) \rangle_{L^2(Q)} \\[1mm]
    && \hspace*{-4cm}
       = \, \langle \partial_t (\widetilde{u} - I_{h_t}^1 P_{h_x}^1
          \widetilde{u}),
          \partial_t (\widetilde{u}_h - I_{h_t}^1 P_{h_x}^1 \widetilde{u})
          \rangle_{L^2(Q)} \\[1mm]
    && + \,
       \langle \nabla_x (\widetilde{u} - I_{h_t}^1 \widetilde{u}),
       \nabla_x \partial_t (\widetilde{u}_h - I_{h_t}^1 P_{h_x}^1
       \widetilde{u}) \rangle_{L^2(Q)} \\[1mm]
    && \hspace*{-4cm}
       = \, \langle \partial_t \widetilde{u} - Q_{h_t}^0 P_{h_x}^1
          \partial_t \widetilde{u} ,
          \partial_t (\widetilde{u}_h - I_{h_t}^1 P_{h_x}^1 \widetilde{u})
          \rangle_{L^2(Q)} \\[1mm]
    && + \, \langle \nabla_x (\widetilde{u} - I_{h_t}^1
       \widetilde{u}),
       \nabla_x \partial_t (\widetilde{u}_h - I_{h_t}^1 P_{h_x}^1
       \widetilde{u}) \rangle_{L^2(Q)} \\[1mm]
    && \hspace*{-4cm} = \,
       \langle \partial_t \widetilde{u} - P_{h_x}^1 \partial_t \widetilde{u},
          \partial_t (\widetilde{u}_h - I_{h_t}^1 P_{h_x}^1 \widetilde{u})
          \rangle_{L^2(Q)} \\[1mm]
    && +
       \langle - \Delta_x (\widetilde{u} - I_{h_t}^1 \widetilde{u}),
       \partial_t (\widetilde{u}_h - I_{h_t}^1 P_{h_x}^1 \widetilde{u})
       \rangle_{L^2(Q)} \\[1mm]
    && \hspace*{-4cm} \leq \, \Big[
           \| \partial_t \widetilde{u} - P_{h_x}^1 \partial_t \widetilde{u}
           \|_{L^2(Q)} +
           \| \Delta_x ( \widetilde{u} - I_{h_t}^1 \widetilde{u})
           \|_{L^2(Q)} \Big]
           \| \partial_t (\widetilde{u}_h - I_{h_t}^1 P_{h_x}^1
           \widetilde{u}) \|_{L^2(Q)},
  \end{eqnarray*}
  i.e.,
  \[
    \| \partial_t (\widetilde{u}_h - I_{h_t}^1 P_{h_x}^1 \widetilde{u})
    \|_{L^2(Q)} \, \leq \,
    \| \partial_t \widetilde{u} - P_{h_x}^1 \partial_t \widetilde{u}
    \|_{L^2(Q)} + \| \Delta_x (\widetilde{u} - I_{h_t}^1
    \widetilde{u}) \|_{L^2(Q)} \, .
  \]
  For $\widetilde{u}_h - I_{h_t}^1 P_{h_x}^1 \widetilde{u} \in X_{0,h}$
  we can use a Friedrich's type inequality to further obtain
  \begin{eqnarray*}
    \| \widetilde{u}_h - I_{h_t}^1 P_{h_x}^1 \widetilde{u} \|_{L^2(Q)}
    & \leq & c \, \| \partial_t (\widetilde{u}_h - I_{h_t}^1 P_{h_x}^1
             \widetilde{u}) \|_{L^2(Q)} \\[1mm]
    & \leq & c \, \Big[ \| \partial_t \widetilde{u} - P_{h_x}^1 \partial_t
             \widetilde{u} \|_{L^2(Q)} +
             \| \Delta_x (\widetilde{u} - I_{h_t}^1 \widetilde{u})
             \|_{L^2(Q)} \Big] \\[1mm]
    & \leq & c \, \Big[ h_x^2 \, | \partial_t \widetilde{u}
             |_{L^2(0,T;H^2(\Omega))} +
             h_t^2 \, |\Delta_x \widetilde{u}|_{H^2(0,T;L^2(\Omega))}
             \Big] .
  \end{eqnarray*}
  Now the assertion follows from
  \[
    \| \widetilde{u} - \widetilde{u}_h \|_{L^2(Q)} \, \leq \,
    \| \widetilde{u} - I_{h_t}^1 P_{h_x}^1 \widetilde{u} \|_{L^2(Q)}
    +
    \| \widetilde{u}_h - I_{h_t}^1 P_{h_x}^1 \widetilde{u} \|_{L^2(Q)},
  \]
  when using, see also the proof of \eqref{Error QI},
  \begin{eqnarray*}
    \| \widetilde{u} - I_{h_t}^1 P_{h_x}^1 \widetilde{u} \|_{L^2(Q)}
    & \leq & \| \widetilde{u} - P_{h_x}^1 \widetilde{u} \|_{L^2(Q)} \, + \,
             \| (\widetilde{u} - I_{h_t}^1 \widetilde{u}) P_{h_x}^1 \|_{L^2(Q)} \\
    & \leq & c_1 \, h_x^2 \, \| \widetilde{u} \|_{L^2(0,T;H^2(\Omega))} \, + \,
             c_2 \, h_t^2 \,
             \| \widetilde{u} \|_{H^2(0,T;L^2(\Omega))} \, .
  \end{eqnarray*}
}

\begin{theorem}\label{Theorem error primal smooth X}
  Let $\widetilde{u}\in X_0$ and $\widetilde{u}_h\in X_{0,h}$ denote the
  unique solutions of \eqref{Heat equation VF hom}
  and \eqref{Heat equation FEM}, respectively. For $\widetilde{u}$
  let all assumptions of Lemma {\rm \ref{Lemma L2 Error primal smooth}}
  be valid. Then there holds the error estimate
  \begin{equation}\label{Error primal X smooth}
    \| \widetilde{u} - \widetilde{u}_h \|_X \,
    \leq \, c(\widetilde{u}) \, h_x,
  \end{equation}
  where $c(\widetilde{u})$ covers all contributions as given in
  \eqref{L2 Error primal smooth}.
\end{theorem}

\myproof{We first note that 
  \[
    \| \widetilde{u} - \widetilde{u}_h \|_X^2 
    \, = \,
    \| \partial_t (\widetilde{u} - \widetilde{u}_h) \|_{Y^*}^2 +
    \| \widetilde{u}-\widetilde{u}_h \|_Y^2
    \, \leq \,
    \| \partial_t (\widetilde{u} - \widetilde{u}_h) \|_{L^2(Q)}^2 +
    \| \widetilde{u}-\widetilde{u}_h \|_Y^2.
  \]
  When using the spatial $H^1_0$ projection $P_{h_x}^1 : H^1_0(\Omega) \to V_{h_x}
  \subset H^1_0(\Omega)$ and the temporal piecewise linear interpolation
  $I_{h_t}^1$ as defined in \eqref{Def Iht}, we obtain 
  \[
    \| \partial_t(\widetilde{u}-\widetilde{u}_h) \|_{L^2(Q)} 
    \, \leq \,
    \| \partial_t(\widetilde{u}-I_{h_t}^1P_{h_x}^1 \widetilde{u}) \|_{L^2(Q)} +
    \| \partial_t(I_{h_t}^1P_{h_x}^1 \widetilde{u})-\widetilde{u}_h \|_{L^2(Q)}, 
  \]
  and as in the proof of Lemma \ref{Lemma L2 Error primal smooth}, we get,
  recall $h_t \sim h_x$,
  \begin{eqnarray*}
    \| \partial_t(I_{h_t}^1P_{h_x}^1 \widetilde{u} - \widetilde{u}_h) \|_{L^2(Q)} 
    & \leq & \| \partial_t \widetilde{u} -
             P_{h_x}^1 \partial_t\widetilde{u} \|_{L^2(Q)} +
             \| \Delta_x(\widetilde{u} - I_{h_t}^1 \widetilde{u}) \|_{L^2(Q)} \\[1mm]
    & \leq & c \, h_x \, \Big( \| \widetilde{u} \|_{H^1(0,T;H^1(\Omega))} +
             \| \widetilde{u} \|_{H^1(0,T;H^2(\Omega))} \Big). 
  \end{eqnarray*}
  Moreover, using
  $\partial_t I_{h_t}^1 \widetilde{u} = Q_{h_t}^0\partial_t \widetilde{u}$,
  and the error estimates \eqref{Error Qht} and \eqref{Error Px L2}, we compute 
  \begin{eqnarray*}
    \| \partial_t (\widetilde{u}-I_{h_t}^1P_{h_x}^1 \widetilde{u}) \|_{L^2(Q)}
    & = & \| \partial_t \widetilde{u} -
          Q_{h_t}^0 P_{h_x}^1 \partial_t\widetilde{u} \|_{L^2(Q)} \\[1mm]
    & \leq & \| \partial_t \widetilde{u} -
             Q_{h_t}^0 \partial_t \widetilde{u} \|_{L^2(Q)} +
             \| Q_{h_t}^0(\partial_t \widetilde{u} -
             P_{h_x}^1 \partial_t \widetilde{u}) \|_{L^2(Q)} \\[1mm]
    & \leq & c \, h_t \, \| \partial_{tt}\widetilde{u} \|_{L^2(Q)} +
             \| \partial_t \widetilde{u} -
             P_{h_x}^1 \partial_t \widetilde{u} \|_{L^2(Q)} \\[1mm]
    & \leq & c \, h_x \, \Big(
             \| \widetilde{u} \|_{H^2(0,T;L^2(\Omega))} +
      \| \widetilde{u} \|_{H^1(0,T;H^1(\Omega))} \Big) .
  \end{eqnarray*}
  Together with the estimate
  \begin{eqnarray*}
    \| \widetilde{u} - \widetilde{u}_h \|_Y
    & \leq & \| \nabla_x (\widetilde{u} - P_{h_x}^1 \widetilde{u})
             \|_{L^2(Q)} + \| \nabla_x (P_{h_x}^1\widetilde{u} -
             \widetilde{u}_h) \|_{L^2(Q)} \\[1mm]
     & \leq & \| \nabla_x (\widetilde{u} - P_{h_x}^1 \widetilde{u})
             \|_{L^2(Q)} + c \, h_x^{-1} \, \| P_{h_x}^1\widetilde{u} -
             \widetilde{u}_h \|_{L^2(Q)} \\[1mm]
     & \leq & \| \nabla_x (\widetilde{u} - P_{h_x}^1 \widetilde{u})
              \|_{L^2(Q)} + c \, h_x^{-1} \, \Big[
              \| P_{h_x}^1\widetilde{u} - \widetilde{u} \|_{L^2(Q)} +
              \| \widetilde{u} - \widetilde{u}_h \|_{L^2(Q)} \Big],
  \end{eqnarray*}
  and using the error estimates \eqref{Error Px H1} and \eqref{Error Px L2}
  gives the assertion \eqref{Error primal X smooth}.}

\noindent
Since all of the above error estimates \eqref{Error primal smooth},
\eqref{L2 Error primal smooth}, and \eqref{Error primal X smooth}
require a sufficiently regular solution $\widetilde{u}$, we now aim to
derive error estimates for less regular solutions. To do so, we
first establish a discrete inf-sup condition in the norm $\| u_h \|_X$
being an equivalent norm to $\nnorm{u_h}_X$ for $u_h \in X_{0,h}$. From
\eqref{upper norm equivalence} we already have that
$\nnorm{u_h}_X \, \leq \, \max\{1,c_Q^{-1}\} \, \| u_h \|_X$ for all
$u_h\in X_{0,h}$. Now we will show the reverse inequality, that is also
outlined in \cite[(5.2.63)]{Andreev:2012}.

\begin{lemma}
  For any $u_h\in X_{0,h}$, there holds 
  \begin{equation}\label{lower norm equivalence}
    \nnorm{u_h}_X \, \geq
    \, c \, \min \{ 1 , h_t^{-1} \, h_x^2 \} \, \| u_h \|_X \, .
  \end{equation}
\end{lemma}

\myproof{When using the Galerkin orthogonality and the error estimate
  \eqref{Error Qht} for the temporal $L^2$ projection $Q_{h_t}^0$, and
  the inverse inequality \eqref{inverse inequality} together with
  a duality argument, we first compute
  \begin{eqnarray*}
    \| u_h \|_Y^2
    & = & \| \nabla_x Q_{h_t}^0 u_h \|_{L^2(Q)}^2 +
          \| \nabla_x (I-Q_{h_t}^0)u_h \|_{L^2(Q)}^2 \\[1mm]
    & \leq & \| Q_{h_t}^0 u_h \|_Y^2 + c \, h_t^2 \,
             \| \partial_t u_h \|_{L^2(0,T;H^1_0(\Omega))}^2 \\[1mm]
    & \leq & \| Q_{h_t}^0 u_h \|_Y^2 + c \,
             \, h_t^2 \, h_x^{-4} \,
             \| \partial_t u_h \|_{L^2(0,T;H^{-1}(\Omega))}^2.
  \end{eqnarray*}
  So, we get
  \begin{eqnarray*}
    \| u_h \|_X^2
    & = & \| \partial_t u_h \|_{Y^*}^2 + \| u_h \|_Y^2 \\[1mm]
    & \leq & \| Q_{h_t}^0 u_h \|_Y^2 + \Big[ 1 +
             c \, h_t^2 \, h_x^{-4} \Big] \,
             \| \partial_t u_h \|_{Y^*}^2 \\[1mm]
    & \leq & \max \{ c_Q^2 , 1 +
             c \, h_t^2 \, h_x^{-4} \} \,
             \Big[ c_Q^{-2} \,
             \| Q_{h_t}^0 u_h \|_Y^2 + \| \partial_t u_h \|_{Y^*}^2 \Big],
  \end{eqnarray*}
  and the assertion follows.}

\noindent
When combining the discrete inf-sup stability condition
\eqref{Heat equation inf sup Error} and the norm equivalence
inequality \eqref{lower norm equivalence}, this results in the
discrete inf-sup stability condition
\begin{equation}\label{eq:discrete inf-sup full X}
  c_{S,h} \, \| u_h \|_X \, \leq \,
  \sup_{0 \neq v_h\in Y_h} \frac{b(u_h,v_h)}{\| v_h \|_Y} \quad
  \mbox{for all} \;  u_h\in X_{0,h},
\end{equation}
where $c_{S,h} := c \, \min \{ 1 , h_t^{-1} h_x^2 \}$. In particular 
we then conclude the Cea type estimate
\begin{equation}\label{Cea primal full X}
  \| \widetilde{u} - \widetilde{u}_h \|_X \, \leq \,
  \frac{\sqrt{2}}{c_{S,h}} \, \inf\limits_{z_h \in X_{0,h}}
  \| \widetilde{u} - z_h \|_X. 
\end{equation}

\begin{remark}\label{Remark error primal X smooth orthotropic}
  Even for smooth solutions $\widetilde{u}$, the quasi best approximation
  error is at best  
  \[
    \| \widetilde{u} - \widetilde{u}_h \|_X \, \leq \,
    \frac{\sqrt{2}}{c_{S,h}} \inf_{z_h\in X_{0,h}}
    \| \widetilde{u}-z_h \|_X \, \leq \,
    \frac{\sqrt{2}}{c_{S,h}} \, h_x \, c(\widetilde{u}). 
  \]
  In the case of orthotropic scaling, i.e., $h_t\simeq h_x$, we have
  $c_{S,h} \simeq h_x$, which would suggest that no convergence is
  observed. However, in \eqref{Error primal X smooth} we derived an estimate
  of order $h_x$ in this setting. This indicates that for $\widetilde{u}$
  satisfying the assumption of Lemma \ref{Lemma L2 Error primal smooth}, the
  discrete inf-sup constant exhibits the behavior
  $c_{S,h} = c_{S,h}(\widetilde{u}) \simeq 1$, i.e., it improves on a certain
  subclass $U$ of functions $\widetilde{u}\in U\subset X_0$, see also
  {\rm \cite{BabushkaNarasimhan:1997}}. 
\end{remark}

\noindent
On the other hand, for parabolic scaling $h_t\simeq h_x^2$, we have
$c_{S,h}\simeq 1$, and with the approximation properties of the space-time
finite element space $X_{0,h}$ we obtain the following result:

\begin{lemma}
  Let $\widetilde{u}\in X_0$ and $\widetilde{u}_h\in X_{0,h}$ denote the
  unique solutions of \eqref{Heat equation VF hom}
  and \eqref{Heat equation FEM}, respectively. Assume
  $\widetilde{u} \in H^{2,1}(Q) \cap H^{3/2}(0,T;H^{-1}(\Omega))$.
  For $h_t \simeq h_x^2$ there holds the error estimate
  \begin{equation}\label{Error primal X H21}
    \| \widetilde{u} - \widetilde{u}_h \|_X \, \leq \, c \, h_x \,  \Big[
    \| \widetilde{u} \|_{H^{2,1}(Q)}^2 +
    \| \widetilde{u} \|_{H^{3/2}(0,T;H^{-1}(\Omega))}^2 \Big]^{1/2} .
  \end{equation}
\end{lemma}

\myproof{We first consider the error estimate \eqref{Error QI 2} for
  $s=2$ and $\sigma=1$,
  \[
    \| \widetilde{u} - Q_{h_x}^1 I_{h_t}^1 \widetilde{u} \|_Y \, \leq \,
    c_1 \, h_x \, \| \widetilde{u} \|_{L^2(0,T;H^2(\Omega))} +
    c_2 \, h_x^{-1} \, h_t \, \| \widetilde{u} \|_{H^1(0,T;L^2(\Omega))},
  \]
  and the error estimate \eqref{Error QI t} for $\tau=0$ and $\varrho=3/2$,
  \[
    \| \widetilde{u} - \widetilde{u}_h \|_{Y^*} \, \leq \,
    c_1 \, h_x \, \| \partial_t \widetilde{u} \|_{L^2(Q)} + c_2 \, h_t^{1/2} \,
    \| \widetilde{u} \|_{H^{3/2}(0,T;H^{-1}(\Omega))} .
  \]
  In particular for $h_t \simeq h_x^2$ and using \eqref{Cea primal full X}
  for $z_h = Q_{h_x}^1 I_{h_t}^1 \widetilde{u}$, the assertion follows.}

\noindent
Similar as for second order elliptic problems, for smoothly bounded or
convex spatial domains $\Omega$, we can now use the
Nitsche trick to conclude an error estimate in $L^2(Q)$.

\begin{lemma}\label{Lemma Nitsche trick}
  Let $\widetilde{u}\in X_0$ and $\widetilde{u}_h\in X_{0,h}$ denote the
  unique solutions of \eqref{Heat equation VF hom}
  and \eqref{Heat equation FEM}, respectively. 
  For $h_t \simeq h_x^2$ there holds the error estimate
  \begin{equation}\label{Error primal L2 H21}
    \| \widetilde{u} - \widetilde{u}_h \|_{L^2(Q)} \, \leq \, c \, h_x \,
    \| \widetilde{u} - \widetilde{u}_h \|_X \, .
  \end{equation}
\end{lemma}

\myproof{Let us consider the adjoint heat equation,
  \begin{equation*}\label{Adjoint Heat equation}
    \begin{array}{rclcl}
      - \partial_t p(x,t) - \Delta_x p(x,t)
      & = & \widetilde{u}(x,t) - \widetilde{u}_h(x,t)
      && \mbox{for} \; (x,t) \in Q, \\[1mm]
      p(x,t)
      & = & 0
      && \mbox{for} \; (x,t) \in \Sigma, \\[1mm]
      p(x,T) & = & 0 && \mbox{for} \; x \in \Omega ,
    \end{array}  
  \end{equation*}
  i.e., $p \in Y$ is the unique solution of the variational formulation
  \[
    b(w,p) \, = \, \langle \widetilde{u} - \widetilde{u}_h , w \rangle_{L^2(Q)}
    \quad \mbox{for all} \; w \in X_0 \, .
  \]
  In particular for $w = \widetilde{u} - \widetilde{u}_h \in X_0$ this gives
  \begin{eqnarray*}
    \| \widetilde{u} - \widetilde{u}_h \|^2_{L^2(Q)}
    & = & \langle \widetilde{u} - \widetilde{u}_h,
          \widetilde{u} - \widetilde{u}_h \rangle_{L^2(Q)} \, = \,
          b(\widetilde{u} - \widetilde{u}_h,p) \\
    & = & b(\widetilde{u} - \widetilde{u}_h,p-Q_{h_x}^1 Q_{h_t}^0 p) \, \leq \,
          \sqrt{2} \, \| \widetilde{u} - \widetilde{u}_h \|_X
          \| p-Q_{h_x}^1 Q_{h_t}^0 p \|_Y \, .
  \end{eqnarray*}
  By parabolic regularity we have that $p\in H^{2,1}(Q)$, since
  $\widetilde{u}-\widetilde{u}_h\in L^2(Q)$, see \eqref{Error primal X H21}.
  Next we use the error estimate \eqref{Error QQ 2} for $s=2$ and $\sigma=1$,
  \[
    \| p-Q_{h_x}^1 Q_{h_t}^0 p \|_Y \, \leq \,
    c_1 \, h_x \, |p|_{L^2(0,T;H^2(\Omega))} + c_2 \, h_x^{-1} \, h_t \,
    \| p \|_{H^1(0,T;L^2(\Omega))} \, \leq \,
    c \, h_x \, \|p\|_{H^{2,1}(Q)},
  \]
  when choosing $h_t \simeq h_x^2$. Now the assertion follows from the
  regularity estimate \eqref{Appendix:Heat f Regularity H21},
  \[
    \|p\|_{H^{2,1}(Q)} \, \leq \,
    \| \widetilde{u} - \widetilde{u}_h \|_{L^2(Q)} \, ,
  \]
and the error estimate \eqref{Error primal X H21}.}

\noindent
Now we are in the position to consider the case when the initial
datum $u_0$ is discontinuous. For this we consider the homogeneous
heat equation \eqref{Heat equation} with $f \equiv 0$, where we can
reformulate the previous error estimates as follows.

\begin{theorem}\label{Theorem disc}
  Let $\widetilde{u}\in X_0$ and $\widetilde{u}_h \in X_{0,h}$ denote the
  unique solutions of \eqref{Heat equation VF hom}
  and \eqref{Heat equation FEM}, respectively. Assume $f \equiv 0$, and
  $u_0 \in \widetilde{H}^s(\Omega) := [L^2(\Omega),H^1_0(\Omega)]_{|s}$
  for some $ s \in [0,1]$. For $h_t \simeq h_x^2$ there hold the error estimates
  \begin{equation}\label{Error primal X Hs}
    \| \widetilde{u} - \widetilde{u}_h \|_X \, \leq \, c \, h_x^s \,
    \| u_0 \|_{\widetilde{H}^s(\Omega)},
  \end{equation}
  and
  \begin{equation}\label{Error primal L2 Hs}
    \| \widetilde{u} - \widetilde{u}_h \|_{L^2(Q)} \, \leq \, c \, h_x^{1+s} \,
    \| u_0 \|_{\widetilde{H}^s(\Omega)}.
  \end{equation}
\end{theorem}

\myproof{For $f \equiv 0$ and assuming $u_0 \in H^1_0(\Omega)$, the error
  estimate \eqref{Error primal X H21} reads
  \[
    \| \widetilde{u} - \widetilde{u}_h \|_X
    \, \leq \, \, h_x \, \Big[ \| \widetilde{u} \|^2_{H^{2,1}(Q)}
    + \| \widetilde{u} \|^2_{H^{3/2}(0,T;H^{-1}(Q))} \Big]^{1/2} \, \leq \,
    c \, h_x \, \| \nabla_x u_0 \|_{L^2(\Omega)},
  \]
  while using \eqref{Cea primal full X} for $z_h = 0$ and the
  regularity estimate \eqref{Appendix:homogeneous heat Regularity L2 H11}
  this gives
  \[
    \| \widetilde{u} - \widetilde{u}_h \|_X \, \leq \,
    c \, \| \widetilde{u} \|_X \, \leq \, c \, \| u_0 \|_{L^2(\Omega)} .
  \]
  Now, \eqref{Error primal X Hs} follows from a space interpolation argument,
  and \eqref{Error primal L2 Hs} follows as in the proof of
  \eqref{Error primal L2 H21}.} 

\noindent
All of the previous error estimates are formulated for the solution
$\widetilde{u}_h$ satisfying the finite element variational formulation
\eqref{Heat equation FEM}, which involves the extension $\widetilde{u}_0$
of the given initial datum $u_0$. Hence, and as in the case of Dirichlet
boundary conditions for elliptic problems, we proceed as follows.
For $\widetilde{u}_0 \in X$ we consider the space-time finite element
approximation $I_{h_t}^1 Q_{h_x}^1 \widetilde{u}_0$, and instead of
\eqref{Heat equation FEM} we consider a perturbed variational formulation
to find $\widetilde{u}_h^* \in X_{0,h}$ such that
\begin{equation}\label{Heat equation FEM pert}
  b(\widetilde{u}_h^* + I_{h_t}^1 Q_{h_x}^1 \widetilde{u}_0,v_h) =
  \langle f , v_h \rangle_Q 
  \quad \mbox{for all} \; v_h \in Y_h ,
\end{equation}
and subtracting \eqref{Heat equation FEM pert} from
\eqref{Heat equation FEM} gives the perturbed Galerkin orthogonality
\begin{equation}\label{pert Galerkin orthogonality}
  b(\widetilde{u}_h - \widetilde{u}_h^*,v_h) =
  b(I_{h_t}^1 Q_{h_x}^1 \widetilde{u}_0-\widetilde{u}_0,v_h) \quad
  \mbox{for all} \; v_h \in Y_h .
\end{equation}
As in the elliptic case we can use the Strang lemma to analyse this
additional error:

\begin{lemma}\label{Lemma primal smooth pert}
  Let $\widetilde{u} \in X_0$ and $\widetilde{u}_h^* \in X_{0,h}$ denote the
  unique solutions of \eqref{Heat equation VF hom} and
  \eqref{Heat equation FEM pert}, respectively. Assume
  $\widetilde{u}, \widetilde{u}_0 \in H^{2,1}(Q) \cap H^2(0,T;H^{-1}(\Omega))$
  such that $\nabla_x \partial_t \widetilde{u}, \nabla_x \partial_t
  \widetilde{u}_0 \in L^2(Q)$ are satisfied, and let $h_t \simeq h_x$.
  Then there holds the error estimate
  \begin{equation}\label{final error primal smooth X}
    \nnorm{\widetilde{u} - \widetilde{u}_h^*}_X \, \leq \,
    c(\widetilde{u},\widetilde{u}_0) \, h_x \, .
  \end{equation}
  If $\widetilde{u}, \widetilde{u}_0 \in
  L^2(0,T;H^2(\Omega)) \cap H^2(0,T;L^2(\Omega))$ and
  $\partial_t \widetilde{u}, \partial_t \widetilde{u}_0 \in
  L^2(0,T;H^2(\Omega))$ as well as $\Delta_x \widetilde{u},
  \Delta_x \widetilde{u}_0\in H^2(0,T;L^2(\Omega))$ is satisfied,
  then
  \begin{equation}\label{final error primal smooth L2}
    \| \widetilde{u} - \widetilde{u}_h^* \|_{L^2(Q)} \, \leq \,
    c(\widetilde{u},\widetilde{u}_0) \, h_x^2 \, ,
  \end{equation}
  and
   \begin{equation}\label{final error primal smooth Y}
    \| \widetilde{u} - \widetilde{u}_h^* \|_Y \, \leq \,
    c(\widetilde{u},\widetilde{u}_0) \, h_x \, .
  \end{equation}
\end{lemma}

\myproof{By the triangle inequality we first have
  \[
    \nnorm{\widetilde{u} - \widetilde{u}_h^*}_X \, \leq \,
    \nnorm{\widetilde{u} - \widetilde{u}_h}_X +
    \nnorm{\widetilde{u}_h - \widetilde{u}_h^*}_X ,
  \]
  and with \eqref{Heat equation inf sup Error} we further obtain
  \begin{eqnarray*}
    \nnorm{\widetilde{u}_h - \widetilde{u}_h^*}_X
    & \leq & \sup\limits_{0 \neq v_h \in Y_h}
             \frac{b(\widetilde{u}_h - \widetilde{u}_h^*,v_h)}{\| v_h \|_Y} \\
     & = & \sup\limits_{0 \neq v_h \in Y_h}
           \frac{b(I_{h_t}^1 Q_{h_x}^1 \widetilde{u}_0-\widetilde{u}_0,v_h)}
           {\| v_h \|_Y} \\
    & \leq & \sqrt{2} \,
             \| \widetilde{u}_0 - I_{h_t}^1 Q_{h_x}^1 \widetilde{u}_0 \|_X
             \, = \, \sqrt{2} \,
             \| \widetilde{u}_0 - Q_{h_x}^1 I_{h_t}^1 \widetilde{u}_0 \|_X,
  \end{eqnarray*}
  since the nodal temporal interpolation $I_{h_t}^1$ and the spatial
  $L^2$ projection $Q_{h_x}^1$ commute. Hence we can use the error estimates
  \eqref{Error QI} and \eqref{Error QI t} to conclude 
  \eqref{final error primal smooth X}. In order to estimate
  $\| \widetilde{u} - \widetilde{u}_h^* \|_{L^2(Q)}$ we can proceed as
  in the proof of Lemma \ref{Lemma L2 Error primal smooth}. But instead
  of \eqref{Proof L2 smooth Galerkin} we have to use the perturbed
  Galerkin orthogonality
  \[
    b(\widetilde{u}_h^*, v_h) \, = \,
    b(\widetilde{u},v_h) +
    b(\widetilde{u}_0-I_{h_t}^1 Q_{h_x}^1 \widetilde{u}_0,v_h)
  \]
  for the particular test function
  $v_h=\partial_t(\widetilde{u}_h-I_{h_t}^1P_{h_x}^1)\widetilde{u} \in Y_h$.
  In any case, the additional term can be estimated similar as in the
  proof of Lemma \ref{Lemma L2 Error primal smooth} to conclude
  \eqref{final error primal smooth L2}. The error estimate
  \eqref{final error primal smooth Y} finally follows as in
  Theorem \ref{Theorem error primal smooth X}.}

\noindent
It remains to consider error estimates for the solution $\widetilde{u}_h^*$
of the perturbed variational formulation \eqref{Heat equation FEM pert},
when the solution $\widetilde{u}$ is less regular, and parabolic
scaling $h_t \simeq h_x^2$ is required.

\begin{lemma}\label{Lemma primal nonsmooth pert}
  Let $\widetilde{u} \in X_0$ and $\widetilde{u}_h^* \in X_{0,h}$ denote the
  unique solutions of \eqref{Heat equation VF hom} and
  \eqref{Heat equation FEM pert}, respectively. Assume
  $\widetilde{u}, \widetilde{u}_0 \in H^{2,1}(Q) \cap
  H^{3/2}(0,T;H^{-1}(\Omega))$. For $h_t \simeq h_x^2$
  there holds 
  \begin{equation}\label{final error primal nonsmooth X}
    \| \widetilde{u} - \widetilde{u}_h^* \|_X \, \leq \,
    \| \widetilde{u} - \widetilde{u}_h \|_X +
    c \, \| \widetilde{u}_0 - I_{h_t}^1 Q_{h_x}^1 \widetilde{u}_0 \|_X,
  \end{equation}
  and
  \begin{equation}\label{final error primal nonsmooth L2}
    \begin{array}{ll}
      &\| \widetilde{u} - \widetilde{u}_h^* \|_{L^2(Q)} \, \leq \,
    c \, h_x \, \Big[ \| \widetilde{u} - \widetilde{u}_h^* \|_X +
    \| \widetilde{u}_0 - I_{h_t}^1 Q_{h_x}^1 \widetilde{u}_0 \|_X \Big]\\ 
    & \hskip 7em 
    \,+\, 
    \| \widetilde{u}_0 - I_{h_t}^1 Q_{h_x}^1 \widetilde{u}_0 \|_{L^2(Q)} 
    \, + \, 
    \|u_0-Q_{h_x}^1 u_0\|_{L^2(\Omega)} .
    \end{array}
  \end{equation}
  
\end{lemma}

\myproof{As in the proof of Lemma \ref{Lemma primal smooth pert} we first
  have
  \[
    \| \widetilde{u} - \widetilde{u}_h^* \|_X \, \leq \,
    \| \widetilde{u} - \widetilde{u}_h \|_X +
    \| \widetilde{u}_h - \widetilde{u}_h^* \|_X,
  \]
  but now we use \eqref{eq:discrete inf-sup full X} for $h_t \simeq h_x^2$
  and \eqref{pert Galerkin orthogonality} to conclude
  \begin{eqnarray*}
    c \, \| \widetilde{u}_h - \widetilde{u}_h^* \|_X
    & \leq & \sup\limits_{0 \neq v_h \in Y_h}
             \frac{b(\widetilde{u}_h - \widetilde{u}_h^*,v_h)}{\| v_h \|_Y} \\
     & = & \sup\limits_{0 \neq v_h \in Y_h}
           \frac{b(I_{h_t}^1 Q_{h_x}^1 \widetilde{u}_0-\widetilde{u}_0,v_h)}
           {\| v_h \|_Y} \, \leq \, \sqrt{2} \,
             \| \widetilde{u}_0 - I_{h_t}^1 Q_{h_x}^1 \widetilde{u}_0 \|_X.
  \end{eqnarray*}
  To prove \eqref{final error primal nonsmooth L2}, we can proceed as
  in the proof of Lemma \ref{Lemma Nitsche trick}, but now we obtain
  \begin{eqnarray*}
    \| \widetilde{u} - \widetilde{u}_h^* \|^2_{L^2(Q)}
    & = & \langle \widetilde{u} - \widetilde{u}_h^*,
          \widetilde{u} - \widetilde{u}_h^* \rangle_{L^2(Q)} \, = \,
          b(\widetilde{u} - \widetilde{u}_h^*,p) \\
    & = & b(\widetilde{u} - \widetilde{u}_h^*,p-Q_{h_x}^1 Q_{h_t}^0 p)
          +
          b(\widetilde{u} - \widetilde{u}_h^*,Q_{h_x}^1 Q_{h_t}^0 p) \\
    & = & b(\widetilde{u} - \widetilde{u}_h^*,p-Q_{h_x}^1 Q_{h_t}^0 p)
          +
          b(I_{h_t}^1 Q_{h_x}^1 \widetilde{u}_0 - \widetilde{u}_0,
          Q_{h_x}^1 Q_{h_t}^0 p),
  \end{eqnarray*}
  and it remains to consider the second term,
  \begin{eqnarray*}
    && b(I_{h_t}^1Q_{h_x}^1\widetilde{u}_0-\widetilde{u}_0,Q_{h_x}^1Q_{h_t}^0p) \\
    && \hspace*{15mm} = \, b(I_{h_t}^1Q_{h_x}^1\widetilde{u}_0-\widetilde{u}_0,
          Q_{h_x}^1Q_{h_t}^0p-p) +
          b(I_{h_t}^1Q_{h_x}^1\widetilde{u}_0-\widetilde{u}_0,p) \\
    && \hspace*{15mm}
       \leq \, \sqrt{2} \, \| \widetilde{u}_0 - I_{h_t}^1 Q_{h_x}^1
       \widetilde{u}_0 \|_X \| p - Q_{h_x}^1 Q_{h_t}^0 p \|_Y +
       b(I_{h_t}^1Q_{h_x}^1\widetilde{u}_0-\widetilde{u}_0,p) ,
  \end{eqnarray*}
  and the first part can be further estimated as in the proof of
  Lemma \ref{Lemma Nitsche trick}. Using integration by parts in space and
  time we obtain
  \begin{eqnarray*}
    b(I_{h_t}^1Q_{h_x}^1\widetilde{u}_0-\widetilde{u}_0,p)
    & = & 
    \langle \partial_t (I_{h_t}^1Q_{h_x}^1\widetilde{u}_0-\widetilde{u}_0),
          p \rangle_Q + \langle \nabla_x
          (I_{h_t}^1Q_{h_x}^1\widetilde{u}_0-\widetilde{u}_0),\nabla_x p
          \rangle_{L^2(Q)} \\
    & = & 
    \langle I_{h_t}^1Q_{h_x}^1\widetilde{u}_0-\widetilde{u}_0 ,
          - \partial_t p - \Delta_x p \rangle_{L^2(Q)} 
          - \langle Q_{h_x}^1 u_0 - u_0,p(0) \rangle_{L^2(\Omega)} \\
    & \leq & 
    \Big[ \| I_{h_t}^1Q_{h_x}^1\widetilde{u}_0-\widetilde{u}_0 \|_{L^2(Q)}
    +
    \| u_0-Q_{h_x}^1 u_0 \|_{L^2(\Omega)} \Big]
             \| p \|_{H^{2,1}(Q)},
  \end{eqnarray*}
  and the assertion follows.}

\noindent
When combining the results of Theorem \ref{Theorem disc}
and Lemma \ref{Lemma primal nonsmooth pert} we finally conclude:
\begin{theorem}\label{Theorem disc pert}
  Let $\widetilde{u}\in X_0$ and $\widetilde{u}_h^* \in X_{0,h}$ denote the
  unique solutions of \eqref{Heat equation VF hom}
  and \eqref{Heat equation FEM pert}, respectively. Assume $f \equiv 0$, and
  $u_0 \in \widetilde{H}^s(\Omega) := [L^2(\Omega),H^1_0(\Omega)]_{|s}$
  for some $ s \in [0,1]$. For $h_t \simeq h_x^2$ there hold the error estimates
  \begin{equation}\label{Error primal X Hs pert}
    \| \widetilde{u} - \widetilde{u}_h^* \|_X \, \leq \, c \, h_x^s \,
    \| u_0 \|_{\widetilde{H}^s(\Omega)},
  \end{equation}
  and
  \begin{equation}\label{Error primal L2 Hs pert}
    \| \widetilde{u} - \widetilde{u}_h^* \|_{L^2(Q)} \, \leq \, c \, h_x^{1+s} \,
    \| u_0 \|_{\widetilde{H}^s(\Omega)}.
  \end{equation}
\end{theorem}

\noindent
By introducing
\[
  u_h := \widetilde{u}_h^* + I_{h_t}^1 Q_{h_x}^1 \widetilde{u}_0 =
  \sum\limits_{i=0}^ {N_t} \sum\limits_{k=1}^{M_x}
  u_k^i \, \psi_i^1(t) \, \varphi_k(x)
\]
we can write the variational formulation \eqref{Heat equation FEM pert} as
\[
  b(u_h,v_h) = \langle f , v_h \rangle_Q \quad \mbox{for all} \; v_h \in Y_h .
\]
For the particular test function
$v_h(x,t) = \psi_i^0(t) \varphi_\ell(x)$, $i=1,\ldots,N_t$,
$\ell=1,\ldots,M_x$, this gives
\begin{eqnarray*}
  && \sum\limits_{k=1}^{M_x}
  \int_{t_{i-1}}^{t_i} \frac{u_k^i-u_k^{i-1}}{h_t} \, dt
  \int_\Omega \varphi_k(x) \, \varphi_\ell(x) \, dx \\
  && \hspace*{5mm} + \sum\limits_{k=1}^{M_x}
     \int_{t_{i-1}}^{t_i} \left( u_k^{i-1} + \frac{t-t_{i-1}}{h_t}
     (u_k^i-u_k^{i-1}) \right) \, dt
     \int_\Omega \nabla_x \varphi_k(x) \cdot \nabla_x \varphi_\ell(x) \, dx \\
  && \hspace*{8cm} = \,
     \int_{t_{i-1}}^{t_i} f(x,t) \varphi_\ell(x) \, dx \, ,
\end{eqnarray*}
which is equivalent to the algebraic system of linear equations,
\begin{equation}
  M_h (\underline{u}^i - \underline{u}^{i-1}) +
  \frac{1}{2} \, h_t \, K_h (\underline{u}^i+\underline{u}^{i-1}) \, = \,
  \underline{f}^i \, , \quad i=1,\ldots,N_t.
\end{equation}
Despite the different evaluation of the right hand side $\underline{f}^i$
this is the Crank--Nicolson scheme, which was introduced in
\cite{CrankNicolson:1947}, to \emph{remove the oscillatory behavior and allow
  for much bigger time-steps} than in a method proposed before by
Richardson \cite{Richardson:1911}. Note that the coefficients
$\underline{u}^0$ are determined from the $L^2$ projection of the initial
datum $u_0$. But it is well known that a piecewise linear continuous
approximation of a discontinuous initial datum $u_0$ produces oscillations
around the discontinuity, see also the discussion in \cite{Osterby2003}.
This motivates the following considerations.

\section{The adjoint space-time FEM for the heat equation}
\label{sec:adjoint discrete}
When considering the variational formulation \eqref{Heat equation VF},
and using integration by parts in time, this gives
\begin{eqnarray*}
  \int_0^T \int_\Omega f(x,t) \, v(x,t) \, dx \, dt
  & = &  \int_0^T \int_\Omega \Big[ \partial_t u(x,t) \, v(x,t) +
        \nabla_x u(x,t) \cdot \nabla_x v(x,t) \Big] \, dx \, dt \\
  && \hspace*{-4.5cm}
     = \, \int_0^T \int_\Omega \Big[ - u(x,t) \, \partial_t v(x,t) +
        \nabla_x u(x,t) \cdot \nabla_x v(x,t) \Big] \, dx \, dt -
        \int_\Omega u_0(x) \, v(x,0) \, dx
\end{eqnarray*}
when using $v(x,T)=0$, $x \in \Omega$. Hence we conclude a
variational formulation to find $u \in Y = L^2(0,T;H^1_0(\Omega))$ such that
\begin{eqnarray}\label{Heat adjoint VF}
  && b_T(u,v) := \int_0^T \int_\Omega \Big[ - u(x,t) \, \partial_t v(x,t) +
     \nabla_x u(x,t) \cdot \nabla_x v(x,t) \Big] \, dx \, dt \\
  && \hspace*{40mm} = \int_0^T \int_\Omega f(x,t) \, v(x,t) \, dx \, dt
     + \int_\Omega u_0(x) \, v(x,0) \, dx \nonumber
\end{eqnarray}
is satisfied for all $v \in X_T := \{ v \in X : v(x,T)=0, x \in \Omega \}$.
Unique solvability of \eqref{Heat adjoint VF} follows as in the
case of the primal variational formulation
\eqref{Heat equation VF}, see also \cite[Lemma 3.1]{LSTY_SISC}.

For the discretization of \eqref{Heat adjoint VF}
we use the ansatz space $Y_h = W_{h_t}^0 \otimes V_{h_x}$ of functions
which are piecewise linear in space but piecewise constant in time, and
the test space $X_{T,h} = \overline{W}_{h_t}^1 \otimes V_{h_x}$, where
$\overline{W}_{h_t}^1 = \mbox{span} \{ \psi_i^1 \}_{i=0}^{N_t-1} \subset
H^1_{,0}(0,T)$ is the space of piecewise linear continuous functions which are
zero at $t=T$. Then we determine $u_h \in Y_h$ such that
\begin{equation}\label{VF adjoint heat}
  b_T(u_h,v_h) = \langle f , v_h \rangle_Q + \langle u_0,
  v_h(0) \rangle_{L^2(\Omega)} \quad \mbox{for all} \; v_h \in X_{T,h}.
\end{equation}
Unique solvability of \eqref{VF adjoint heat} 
is a direct consequence of the discrete inf-sup condition
\eqref{eq:discrete inf-sup full X}.

\begin{lemma}
  For $Y_h = W_{h_t}^0 \otimes V_{h_x}$ and
  $X_{T,h} = \overline{W}_{h_t}^1 \otimes V_{h_x}$ there holds the discrete
  inf-sup stability condition
  \begin{equation}\label{adjoint heat inf-sup sharp}
    c_{S,h} \,
    \| u_h \|_Y \, \leq \, \sup\limits_{0 \neq v_h \in X_{T,h}}
    \frac{b_T(u_h,v_h)}{\| v_h \|_X} \quad
    \mbox{for all} \; u_h \in Y_h,
  \end{equation}
  where $c_{S,h} = c \, \min \{ 1 , h_t^{-1} \, h_x^2 \}$.
\end{lemma}

\myproof{Consider the time reversal operator $\iota_T v(t) = v(T-t),$ which
  is an isometry. For any given $u_h \in Y_h$ we first compute
  $\overline{v}_h \in X_{0,h} = \iota_T(X_{T,h})$ as unique solution
  of the variational formulation
  \[
    b(\overline{v}_h,\iota_T w_h) 
    \, = \, 
    \langle \nabla_x \iota_T u_h ,
    \nabla_x \iota_Tw_h \rangle_{L^2(Q)} \quad \mbox{for all} \; w_h \in Y_h ,
  \]
  and using the discrete inf-sup stability condition
  \eqref{eq:discrete inf-sup full X} this gives
  \begin{eqnarray*}
    c \, \min \{ 1 , h_t^{-1} \, h_x^2 \} \, \| \overline{v}_h \|_X
    & \leq & \sup\limits_{0 \neq w_h \in Y_h}
             \frac{b(\overline{v}_h,\iota_T w_h)}{\|w_h \|_Y} \\
    & = & \sup\limits_{0 \neq w_h \in Y_h}
          \frac{\langle \nabla_x \iota_T u_h , \nabla_x \iota_Tw_h
          \rangle_{L^2(Q)}}{\|w_h \|_Y} \, \leq \, \| \iota_Tu_h \|_Y 
          = \|u_h\|_Y.
  \end{eqnarray*}
  Now let $\widehat{v}_h(t) := \iota_T\overline{v}_h(t)$, $t \in (0,T)$,
  then we can write, see Lemma \ref{lem:discrete inf-sup equality},
  \begin{eqnarray*}
    b_T(u_h,\widehat{v}_h)
    & = & b(\overline{v}_h,\iota_T u_h) \, = \,
          \langle \nabla_x \iota_Tu_h , \nabla_x \iota_T u_h \rangle_{L^2(Q)} \\
    & = & \| u_h \|_Y^2 \, \geq \,
          c \, \min \{ 1 , h_t^{-1} \, h_x^2 \} \, \| u_h \|_Y
          \| \overline{v}_h \|_X,
  \end{eqnarray*}
  implying the assertion, noting that
  $\| \overline{v}_h \|_X = \| \widehat{v} \|_X$.}

\begin{remark}
  In fact, one can show that the discrete inf-sup constants for the primal and
  the adjoint problem are equal, see Lemma \ref{lem:discrete inf-sup equality},
  i.e., 
  \begin{equation*}
    \inf_{u_h\in X_{0,h}} \sup_{0\neq v_h\in Y_h}
    \frac{b(u_h,v_h)}{\| u_h\|_X\|v_h\|_Y} =
    \inf_{v_h\in Y_h} \sup_{0\neq u_h\in X_{T,h}}
    \frac{b_T(v_h,u_h)}{\|v_h\|_Y \|u_h\|_X}.   
  \end{equation*}
  Thus, we assume if $u$ satisfies the assumption of
  Lemma {\rm \ref{Lemma L2 Error primal smooth}}, we again have in the
  case of orthotropic scaling $h_t\simeq h_x$ that $c_{S,h}\simeq 1$. 
\end{remark}

\noindent
As a direct consequence of the discrete inf-sup stability condition
\eqref{adjoint heat inf-sup sharp}, and as in \eqref{Cea primal full X},
we conclude Cea's lemma
\begin{equation}\label{adjoint Cea}
  \| u - u_h \|_Y \, \leq \,
  \frac{\sqrt{2}}{c_{S,h}} \inf\limits_{w_h \in Y_h} \| u - w_h \|_Y \, .
\end{equation}
When combining this with the error estimate \eqref{Error QQ 2}
for $s=2$ and $\sigma=1$, we conclude the following result:

\begin{theorem}\label{thm:Error adjoint Y smooth}
  Let $u \in Y$ and $u_h \in Y_h$ be the unique solutions of the variational
  formulations \eqref{Heat adjoint VF} and \eqref{VF adjoint heat},
  respectively. Assume $u \in H^{2,1}(Q)$, and let $h_t \simeq h_x^2$,
  or let $u$ satisfy all the assumptions of
  Lemma {\rm \ref{Lemma L2 Error primal smooth}}.
  Then there holds the error estimate
  \begin{equation}
    \| u - u_h \|_Y \, \leq \, c \, h_x \, |u|_{H^{2,1}(Q)} .
  \end{equation}
\end{theorem}

\noindent
It remains to derive an estimate for the error in $L^2(Q)$:

\begin{lemma}\label{lem:error estimate L^2 adjoint}
  For $u \in Y$ and $u_h \in Y_h$ being the unique solutions of the variational
  problems \eqref{Heat adjoint VF} and \eqref{VF adjoint heat}, there
  holds the error estimate
  \begin{equation}\label{Error adjoint L2}
    \| u - u_h \|_{L^2(Q)} \, \leq \, c \, \Big[
    \| u - P_{h_x}^1 u \|_{L^2(Q)} + \| u - Q_{h_t}^0 u \|_{L^2(Q)} +
    h_x^{-1} \, h_t \, \| \nabla_x (u-Q_{h_t}^0u) \|_{L^2(Q)} \Big] .
  \end{equation}
\end{lemma}

\myproof{From the variational formulations \eqref{Heat adjoint VF} and
  \eqref{VF adjoint heat} we first conclude the Galerkin orthogonality
  \[
    - \langle u_h , \partial_t v_h \rangle_Q +
    \langle \nabla_x u_h , \nabla_x v_h \rangle_{L^2(Q)}
    \, = \,
    - \langle u , \partial_t v_h \rangle_Q +
    \langle \nabla_x u , \nabla_x v_h \rangle_{L^2(Q)} , \quad
    v_h \in X_{T,h}.
  \]
  For $w_h \in Y_h$, define
  \[
    \overline{v}_h(x,t) := \int_t^T w_h(x,s) \, ds, \quad
    \overline{v}_h \in X_{T,h}, \quad \partial_t \overline{v}_h(x,t) =
    - w_h(x,t) ,
  \]
  where we have
  \begin{eqnarray*}
    && \int_0^T \int_\Omega \nabla_x w_h(x,t) \cdot
       \nabla_x \overline{v}_h(x,t) \, dx \, dt
       \, = \, - \int_0^T \int_\Omega \nabla_x \partial_t
       \overline{v}_h(x,t) \cdot \nabla_x \overline{v}_h(x,t) \, dx \, dt \\
    && \hspace*{2cm} = \, - \frac{1}{2} \int_0^T \frac{d}{dt} \int_\Omega
       |\nabla_x \overline{v}_h(x,t)|^2 \, dx \, dt \\
    && \hspace*{2cm} = \,
       \left. - \frac{1}{2} \, \| \nabla_x \overline{v}_h(t)
          \|^2_{L^2(\Omega)} \right|_0^T \, = \,
          \frac{1}{2} \, \| \nabla_x \overline{v}_h(0) \|^2_{L^2(\Omega)} 
          \, \geq \, 0.
  \end{eqnarray*}
  In particular for $w_h = u_h - P_{h_x}^1 Q_{h_t}^0u \in Y_h$ we can write
  \begin{eqnarray*}
    \| u_h - P_{h_x}^1 Q_{h_t}^0 u \|^2_{L^2(Q)}
    & = & \langle u_h - P_{h_x}^1 Q_{h_t}^0 u ,
          u_h - P_{h_x}^1 Q_{h_t}^0 u \rangle_{L^2(Q)} \\[1mm]
    && \hspace*{-2cm}
       = \, - \langle u_h - P_{h_x}^1 Q_{h_t}^0 u , \partial_t \overline{v}_h
          \rangle_{L^2(Q)} \\[1mm]
    && \hspace*{-2cm} \leq
       \, - \langle u_h - P_{h_x}^1 Q_{h_t}^0 u , \partial_t \overline{v}_h
           \rangle_{L^2(Q)} + \langle \nabla_x(u_h-P_{h_x}^1 Q_{h_t}^0u),
           \nabla_x \overline{v}_h \rangle_{L^2(Q)} \\[1mm]
    && \hspace*{-2cm}
       = \, - \langle u - P_{h_x}^1 Q_{h_t}^0 u , \partial_t \overline{v}_h
           \rangle_{L^2(Q)} + \langle \nabla_x(u-P_{h_x}^1 Q_{h_t}^0u),
           \nabla_x \overline{v}_h \rangle_{L^2(Q)} .
  \end{eqnarray*}
  For the first term we have
  \begin{eqnarray*}
    - \langle u - P_{h_x}^1 Q_{h_t}^0 u , \partial_t \overline{v}_h
    \rangle_{L^2(Q)}
    & = & - \langle u - P_{h_x}^1 u , \partial_t \overline{v}_h
          \rangle_{L^2(Q)} \\[1mm]
    & \leq & \| u - P_{h_x}^1 u \|_{L^2(Q)}
           \| \partial_t \overline{v}_h \|_{L^2(Q)} \\[1mm]
    & = & \| u - P_{h_x}^1 u \|_{L^2(Q)} \| u_h - P_{h_x}^1 Q_{h_t}^0 u
          \|_{L^2(Q)},
  \end{eqnarray*}
  and for the second term we have
  \begin{eqnarray*}
    \langle \nabla_x(u-P_{h_x}^1 Q_{h_t}^0u),
    \nabla_x \overline{v}_h \rangle_{L^2(Q)}
    & = & \langle \nabla_x(u-Q_{h_t}^0u), \nabla_x \overline{v}_h 
          \rangle_{L^2(Q)} \\[1mm]
    && \hspace*{-2cm} = \, \langle \nabla_x(u-Q_{h_t}^0u),
          \nabla_x (\overline{v}_h-Q_{h_t}^0 \overline{v}_h)      
          \rangle_{L^2(Q)} \\[1mm]
    && \hspace*{-2cm} \leq \, \| \nabla_x(u-Q_{h_t}^0u) \|_{L^2(Q)}
             \| \nabla_x (\overline{v}_h-Q_{h_t}^0 \overline{v}_h)         
             \|_{L^2(Q)} \\[1mm]
    && \hspace*{-2cm} \leq \, c \, h_x^{-1} \, \| \nabla_x(u-Q_{h_t}^0u) \|_{L^2(Q)}
             \| \overline{v}_h-Q_{h_t}^0 \overline{v}_h \|_{L^2(Q)} \\[1mm]
    && \hspace*{-2cm}
       \leq \, c \, h_x^{-1} \, h_t \, \| \nabla_x(u-Q_{h_t}^0u) \|_{L^2(Q)}
             \| \partial_t \overline{v}_h \|_{L^2(Q)} \\[1mm]
    && \hspace*{-2cm}
       = \, c \, h_x^{-1} \, h_t \, \| \nabla_x(u-Q_{h_t}^0u) \|_{L^2(Q)}
           \| u_h - P_{h_x}^1 Q_{h_t}^0 u \|_{L^2(Q)} .
  \end{eqnarray*}
  Hence we conclude
  \[
    \| u_h - P_{h_x}^1 Q_{h_t}^0 u \|_{L^2(Q)} \, \leq \,
    \| u - P_{h_x}^1 u \|_{L^2(Q)} +
    c \, h_x^{-1} \, h_t \, \| \nabla_x(u-Q_{h_t}^0u) \|_{L^2(Q)} .
  \]
  With
  \[
    \| u - u_h \|_{L^2(Q)} \, \leq \,
    \| u - P_{h_x}^1 Q_{h_t}^0 u \|_{L^2(Q)} +
    \| u_h - P_{h_x}^1 Q_{h_t}^0 u \|_{L^2(Q)} 
  \]
  and
  \begin{eqnarray*}
    \| u - P_{h_x}^1 Q_{h_t}^0 u \|_{L^2(Q)}
    & \leq & \| u - Q_{h_t}^0 u \|_{L^2(Q)} +
             \| Q_{h_t}^0 (u-P_{h_x}^1u) \|_{L^2(Q)} \\
    & \leq & \| u - Q_{h_t}^0 u \|_{L^2(Q)} +
             \| u-P_{h_x}^1u \|_{L^2(Q)},
  \end{eqnarray*}
  the assertion follows.}

\begin{corollary}\label{cor:Error L2 adjoint smooth}
  Assume $u \in H^{2,1}(Q)$ such that $\nabla_x \partial_t u \in L^2(Q)$
  is satisfied. Then there holds the error estimate
  \begin{equation}
    \| u - u_h \|_{L^2(Q)} \, \leq \, c \, \Big[
    h_x^2 \, \| u \|_{L^2(0,T;H^2(\Omega))}
    +
    h_t \, \| u \|_{H^1(0,T;L^2(\Omega))} +
    h_x^{-1} \, h_t^2 \, \| \nabla_x \partial_t u  \|_{L^2(Q)} \Big].
  \end{equation}
  In particular for $h_t \simeq h_x$ this gives
  \begin{equation}
    \| u - u_h \|_{L^2(Q)} \, \leq \, c \, h_x \, \Big[
    h_x \, \| u \|_{L^2(0,T;H^2(\Omega))}
    +
    \| u \|_{H^1(0,T;L^2(\Omega))} +
    \| \nabla_x \partial_t u  \|_{L^2(Q)} \Big],
  \end{equation}
  while for $h_t \simeq h_x^2$ we obtain
  \begin{equation}
    \| u - u_h \|_{L^2(Q)} \, \leq \, c \, h_x^2 \, \Big[
    \| u \|_{L^2(0,T;H^2(\Omega))}
    +
    \| u \|_{H^1(0,T;L^2(\Omega))} +
    h_x \, \| \nabla_x \partial_t u  \|_{L^2(Q)} \Big].
  \end{equation}
  Even for a smooth solution $u$ we therefore need to use the parabolic
  scaling $h_t \simeq h_x^2$ to ensure second order convergence in $L^2(Q)$ when
  using the adjoint formulation.
\end{corollary}

\noindent
The need of a parabolic scaling to get higher order rates in $L^2(Q)$, also
stems from the fact that for the adjoint formulation the ansatz space $Y_h$
consists of only piecewise constant functions in time which do not give higher order
convergence, even for smooth functions. This can be circumvented using a
reconstruction operator in a post-processing step. To obtain a piecewise
linear, globally continuous solution in time, we consider the following
reconstruction operator $\mathcal{R}_{h_t}: W_{h_t}^0 \to W_{h_t}^1$ where
for each $u_{h_t}\in W_{h_t}^0$ given as 
\[
  u_{h_t}(t) = \sum_{\ell=1}^N u_\ell \psi_\ell^0(t)
\]
the reconstruction is defined as
\begin{equation}\label{eq:time reconstruction}
  (\mathcal{R}_{h_t} u_{h_t})(t) = \sum_{k=1}^N \bar u_k\psi^1_k(t),           
\end{equation}
where 
\[
  \bar u_k = \frac{1}{2} \, (u_k + u_{k+1}) \quad \mbox{for} \; k=1,\ldots,N-1,
  \quad
  \bar u_N = \frac{1}{2} \, (3 u_N - u_{N-1}).
\]
For ease of presentation, we assume that the temporal mesh is uniform. In the
non-uniform case, define 
\begin{eqnarray*}
  \bar u_k
  & = & \frac{h_{t,k} \, u_k + h_{t,k+1} \, u_{k+1}}{h_{t,k+1} + h_{t,k}}
        \quad \mbox{for} \; k=1,\ldots,N-1,\\
  \bar u_N
  & = & 
        \frac{(2h_{h,N} + h_{t,N-1}) \, u_N - h_{t,N} \, u_{N-1}}
        {h_{t,N-1} + h_{t,N}} \, .
\end{eqnarray*}
By definition we set $({\mathcal{R}}_{h_t}u_{h_t})(0) = 0$, i.e., $\bar u_0 = 0$.
Similar reconstructions are well known for DG in time schemes, see, e.g.,
\cite{Gomez:2026}, for an excellent overview. 
The reconstruction operator satisfies the following properties. 

\begin{lemma}\label{lem:reconstruction properties}
  For any $u_{h_t}\in W_{h_t}^0$ it holds that 
  \begin{equation}\label{eq:stability of the reconstruction}
    \| \mathcal{R}_{h_t}u_{h_t} \|_{L^2(t_{\ell-1},t_{\ell})} \, \leq \,
    c \, \| u_{h_t} \|_{L^2(\omega_\ell)} \quad  \mbox{for all} \; \ell=1,\ldots,N,
  \end{equation}
  where the element patch $\omega_\ell$ is defined as 
  \[
    \omega_\ell = \begin{cases}
      (t_0,t_2),& \ell=1,\\
      (t_{\ell-2},t_{\ell+1}),& \ell=2,\ldots,N-1,\\
      (t_{N-2},t_N),& \ell=N.
    \end{cases}
  \]
  Moreover, for any $u\in H^2(0,T)$ satisfying $u(0)=0$ we have 
  \begin{equation}\label{eq:error estimate reconstruction}
    \| u - \mathcal{R}_{h_t} Q_{h_t}^0 u \|_{L^2(0,T)}
    \, \leq \, c \, h_t^2 \, |u|_{H^2(0,T)}. 
  \end{equation}
\end{lemma}

\myproof{For a piecewise constant function $u_{h_t}\in W_{h_t}^0$, given as 
  \[
    u_{h_t}(t) = \sum_{\ell=1}^{N} u_\ell \psi_{\ell}^0(t),
  \]
  we have  
  \[
    \| u_{h_t} \|_{L^2(t_{\ell-1},t_\ell)}^2 \, = \, h_t \, u_\ell^2,
    \quad \ell=1,\ldots,N. 
  \]
  The reconstruction ${\mathcal{R}}_{h_t}u_{h_t}$ is a piecewise linear,
  globally continuous function,
  \[
    (\mathcal{R}_{h_t} u_{h_t})(t) \, = \, \sum_{k=1}^N \bar u_k \psi^1_k(t).
  \]
  On each temporal element $(t_{\ell-1},t_\ell)$, we can compute 
  \begin{eqnarray*}
    \| \mathcal{R}_{h_t}u_{h_t} \|_{L^2(t_{\ell-1},t_\ell)}^2
    & = & \int_{t_{\ell-1}}^{t_\ell} \left[
          \frac{t_\ell-t}{t_\ell-t_{\ell-1}} \, \bar u_{\ell-1} +
          \frac{t-t_{\ell-1}}{t_\ell - t_{\ell-1}} \, \bar u_\ell \right]^2 \,
          dt \\
    & = & \frac{1}{3} \, h_t \, \big(\bar u_{\ell-1}^2 +
          \bar u_{\ell-1}\bar u_\ell + \bar u_\ell^2\big) \, \leq \,
          \frac{1}{2} \, h_t \, \Big( \bar{u}_{\ell-1}^2 + \bar{u}_\ell^2 \Big) .
  \end{eqnarray*}
  For all interior nodes, i.e., for $\ell=1,\ldots,N-1$, we have 
  \[
    \bar u_\ell^2 \, = \,
    \left( \frac{1}{2} \, (u_{\ell} + u_{\ell+1}) \right)^2 \, \leq \,
    \frac{1}{2} \, \big( u_\ell^2 + u_{\ell+1}^2 \big),
  \]
  while for $\ell=N$ we get 
  \begin{eqnarray*}
    \bar u_N^2 
    & = & \left(\frac{1}{2} \, (3 u_N - u_{N-1}) \right)^2 \\
    & = & \frac{1}{4} \,
          \left(
          \left(
          \begin{array}{cc}
            9 & -3 \\
            -3 & 1
          \end{array}
          \right)
                 \left( \begin{array}{c} u_N \\ u_{N-1} \end{array} \right),
    \left( \begin{array}{c} u_N \\ u_{N-1} \end{array} \right)
    \right)
    \, \leq \, \frac{5}{2} \, \big( u_N^2 + u_{N-1}^2 \big) .
  \end{eqnarray*}
  Hence we obtain
  \[
    \| \mathcal{R}_{h_t} u_{h_t} \|_{L^2(t_{0},t_1)}^2 \, \leq \,
    \frac{1}{2} \, h_t \, \bar{u}_1^2
    \, \leq \, \frac{1}{4} \, h_t \, \big( u_1^2 + u_2^2 \big)
    \, = \, \frac{1}{4} \, \| u_{h_t} \|_{L^2(t_0,t_2)}^2 ,
  \]
  and for $\ell=2,\ldots,N-1$, we have 
  \[
    \| \mathcal{R}_{h_t} u_{h_t} \|_{L^2(t_{\ell-1},t_\ell)}^2
    \leq \frac{1}{2} \, h_t \,
    \big( \bar{u}_{\ell-1}^2 + \bar{u}_\ell^2 \big)
    \leq \frac{1}{4} \, h_t \,
    \big( u_{\ell-1}^2 + 2 u_\ell^2 + u_{\ell+1}^2 \big) 
    \leq \frac{1}{2} \, \| u_{h_t} \|_{L^2(t_{\ell-2},t_{\ell+1})}^2 .  
  \]
  It remains to consider $\ell=N$,  
  \[
    \| \mathcal{R}_{h_t}u_{h_t} \|_{L^2(t_{N-1},t_N)}^2
    \leq \frac{1}{2} \, h_t \, \big( \bar{u}_{N-1}^2 + \bar{u}_N^2 \big)
    \leq
    \frac{3}{2} \, h_t \, \big( u_{N-1}^2 + u_N^2 \big)
     \leq \frac{3}{2} \, \| u_{h_t} \|_{L^2(t_{N-2},t_N)}^2. 
   \]
   This concludes the proof of the stability estimate
   \eqref{eq:stability of the reconstruction}.

   In order to prove \eqref{eq:error estimate reconstruction} we consider
   an arbitrary linear function $q(t) = a+b \, t$ for which we compute 
   \[
     (Q_{h_t}^0q)(t) = \sum\limits_{\ell=1}^N q_\ell \psi_\ell^0(t), \quad
     q_\ell \, = \,
     \frac{1}{h_t} \int_{t_{\ell-1}}^{t_\ell} q(t) \, dt
     \, = \,
     a + \frac{1}{2} \, b \, (t_\ell + t_{\ell-1}) .
   \]
   With this we further obtain
   \[
     ({\mathcal{R}}_{h_t} Q_{h_t}^0 q)(t) =
     \sum\limits_{\ell=1}^N \bar{q}_\ell \psi_\ell^1(t),
   \]
   with, for $\ell=1,\ldots,N-1$,
   \[
     \bar{q}_\ell = \frac{1}{2} \, (q_\ell + q_{\ell+1}) =
     \frac{1}{2} \left( a + \frac{1}{2} \, b \, (t_\ell + t_{\ell-1})
       + a + \frac{1}{2} \, b \, (t_{\ell+1} + t_\ell) \right) =
     a + b \, t_\ell = q(t_\ell),
   \]
   where we have used $t_{\ell \pm 1} = t_\ell \pm h_t$. Moreover,
   for $\ell=N$ we conclude
   \begin{eqnarray*}
     \bar{q}_N
     & = & \frac{1}{2} \, (3q_N-q_{N-1}) =
           \frac{1}{2} \left( 3 \left( a + \frac{1}{2} \, b \,
           (t_N+t_{N-1}) \right) - \left( a + \frac{1}{2} \, b \,
           (t_{N-1}+t_{N-2}) \right) \right) \\
     & = & a + b \, t_N \, = \, q(t_N) .
   \end{eqnarray*}
   Hence, on each interval $(t_{\ell-1},t_\ell)$, $\ell=2,\ldots,N$, it
   holds that 
   \[
     (\mathcal{R}_{h_t} Q_{h_t}^0 q)(t) = q(t)
   \]
   for any affine linear function $q(t) = a+b \, t$, while on the first
   interval, $(t_0,t_1)$, due to the homogeneous initial condition, we get  
   \[
     (\mathcal{R}_{h_t}Q_{h_t}^0 q)(t) = q(t)
     \quad \mbox{for all} \; q(t) = b \, t.
   \]
   Therefore we get for any $u \in H^2(\omega_\ell)$, using the stability
   estimate \eqref{eq:stability of the reconstruction}, that  
   \begin{eqnarray*}
     \| u - \mathcal{R}_{h_t}Q_{h_t}^0 u \|_{L^2(t_{\ell-1},t_\ell)}
     & = & \| u - q + \mathcal{R}_{h_t}Q_{h_t}^0 (q-u)
           \|_{L^2(t_{\ell-1},t_\ell)}\\
     & & \hspace*{-4cm} \leq \,
         \| u-q \|_{L^2(t_{\ell-1},t_\ell)} +
         \| \mathcal{R}_{h_t}Q_{h_t}^0 (q-u) \|_{L^2(t_{\ell-1,t_\ell})} \,
         \leq \, c \, \| u- q \|_{L^2(\omega_\ell)}. 
   \end{eqnarray*}
   By standard approximation results, since $q$ was arbitrary,
   we get for $\ell=1$  
   \[
     \| u - \mathcal{R}_{h_t}Q_{h_t}^0 u \|_{L^2(t_0,t_1)} 
     \, \leq \, c \, \inf_{q(t) = b \, t} \| u-q \|_{L^2(0,t_2)} 
     \, \leq \, c \, (2h_t)^2 \, |u|_{H^2(\omega_\ell)}, 
   \]
   and for $\ell=2,\ldots,N$ 
   \[
     \| u-\mathcal{R}_{h_t}Q_{h_t}^0 u \|_{L^2(t_{\ell-1},t_\ell)}
     \, \leq \, c \, \inf_{q(t) = a+b\,t } \| u-q \|_{L^2(\omega_\ell)} 
     \, \leq c \, (3h_t)^2 \, |u|_{H^2(\omega_\ell)}. 
   \]
   Summing over all $\ell=1,\ldots,N$ completes the proof.}

 \begin{theorem}\label{thm:error reconstruction in time}
   Let $u\in Y$ and $u_h\in Y_h$ be the unique solutions of
   \eqref{Heat adjoint VF} and \eqref{VF adjoint heat}, respectively.
   If $u \in H^2(0,T;L^2(\Omega))\cap L^2(0,T;H^2(\Omega))$,
   $\nabla_x u\in H^2(0,T;L^2(\Omega))$ and
   $\partial_t \Delta_x u\in L^2(Q)$ are satisfied, then for
   $h_t\simeq h_x$ we have the error estimate
   \[
     \| u - \mathcal{R}_{h_t} u_h \|_{L^2(Q)} \leq c(u) \, h_x^2,
   \]
   where $c(u)$ is a constant depending on $u$.
 \end{theorem}

 \myproof{We split the error 
   \[
     \| u - \mathcal{R}_{h_t}u_h \|_{L^2(Q)} \, \leq \,
     \| u - \mathcal{R}_{h_t}Q_{h_t}^0 P_{h_x}^1 u \|_{L^2(Q)} +
     \| \mathcal{R}_{h_t} Q_{h_t}^0 P_{h_x}^1 u - \mathcal{R}_{h_t} u_h
     \|_{L^2(Q)}. 
   \]
   For the first term, using a Poincar\'{e} type inequality, the
   stability of $P_{h_x}^1$, as well as the error estimates
   \eqref{Error Px L2} and \eqref{eq:error estimate reconstruction}, we get 
   \begin{eqnarray*}
     \| u - \mathcal{R}_{h_t} Q_{h_t}^0 P_{h_x}^1 u \|_{L^2(Q)} 
     & \leq & \| u - P_{h_x}^1 u \|_{L^2(Q)} +
       \| P_{h_x}^1(\mathcal{R}_{h_t}Q_{h_t}^0 u-u) \|_{L^2(Q)} \\
     & \leq & \| u-P_{h_x}^1 u \|_{L^2(Q)} +
              c \, \| \nabla_x P_{h_x}^1 (\mathcal{R}_{h_t}Q_{h_t}^0 u-u)
              \|_{L^2(Q)} \\
     & \leq & \| u-P_{h_x}^1 u \|_{L^2(Q)} +
              c \, \| \nabla_x (\mathcal{R}_{h_t}Q_{h_t}^0 u-u) \|_{L^2(Q)} \\
     & \leq & c \, \Big(
              h_x^2 \, |u|_{L^2(0,T;H^2(\Omega))} +
              h_t^2 \, |\nabla_x u|_{H^2(0,T;L^2(\Omega))} \Big) . 
   \end{eqnarray*}
   For the second term, using the stability estimate
   \eqref{eq:stability of the reconstruction}, we have that 
   \[
     \| \mathcal{R}_{h_t}Q_{h_t}^0 P_{h_x}^1 u -
     \mathcal{R}_{h_t}u_h \|_{L^2(Q)} \, \leq \, c \,
     \| Q_{h_t}^0 P_{h_x}^1 u - u_h \|_{L^2(Q)}. 
   \]
   Now, we proceed similar as in the proof of
   Lemma \ref{lem:error estimate L^2 adjoint}. When
   using the Galerkin orthogonality
   \[
     b_T(u,v_h) = b_T(u_h,v_h) \quad \mbox{for all} \, v_h\in X_{T,h},
   \]
   we compute that 
   \begin{eqnarray*}
     && \hspace*{-3mm} b_T(Q_{h_t}^0 P_{h_x}^1u-u_h,v_h) 
            \, = \, b_T(Q_{h_t}^0P_{h_x}^1u -u,v_h) \\[1mm]
     && = \, - \langle Q_{h_t}^0 P_{h_x}^1 u - u,
        \partial_t v_h \rangle_{L^2(Q)} +
        \langle \nabla_x (Q_{h_t}^0P_{h_x}^1u - u),\nabla_x v_h
        \rangle_{L^2(Q)} \\[1mm]
     && = \, - \langle Q_{h_t}^0 P_{h_x}^1u - P_{h_x}^1u,\partial_t v_h
        \rangle_{L^2(Q)} - \langle P_{h_x}^1u-u,\partial_t v_h \rangle_{L^2(Q)}
     \\
     && \hspace*{5cm}
        + \, \langle \nabla_x (Q_{h_t}^0u - u),\nabla_x v_h \rangle_{L^2(Q)}\\[1mm]
     && = \, - \langle P_{h_x}^1 u - u,\partial_t v_h \rangle_{L^2(Q)} +
        \langle \nabla_x (Q_{h_t}^0 u-u),\nabla_x (v_h-Q_{h_t}^0 v_h)
        \rangle_{L^2(Q)} 
   \end{eqnarray*}
   In particular for
   \[
     \overline{v}_h(x,t) := \int_t^T [Q_{h_t}^0 P_{h_x}^1
     u(x,s)-u_h(x,s)] \, ds, \quad \overline{v}_h \in X_{T,h},
   \]
   we get $-\partial_t \overline{v}_h = Q_{h_t}^0 P_{h_x}^1 u-u_h$, and     
   \[
     - \langle Q_{h_t}^0 P_{h_x}^1 u - u_h,\partial_t \overline{v}_h
     \rangle_{L^2(Q)} \, = \, \| Q_{h_t}^0 P_{h_x}^1 u-u_h \|_{L^2(Q)}^2,  
   \]
   as well as 
   \begin{eqnarray*}
     \langle \nabla_x(Q_{h_t}^0 P_{h_x}^1 u-u_h),\nabla_x \overline{v}_h
     \rangle_{L^2(Q)} 
     & = & - \langle \nabla_x \partial_t \overline{v}_h,
           \nabla_x \overline{v}_h \rangle_{L^2(Q)} \\
     && \hspace*{-3cm} = \, - \frac{1}{2} \int_0^T \frac{d}{dt}
           \| \nabla_x \overline{v}_h(t) \|_{L^2(\Omega)}^2 \, dt \, = \,
           \frac{1}{2} \, \| \nabla_x \overline{v}_h(0) \|_{L^2(\Omega)}^2
           \, \geq \, 0. 
   \end{eqnarray*}
   Hence, we get 
   \begin{eqnarray*}
     && \| Q_{h_t}^0 P_{h_x}^1 u-u_h \|_{L^2(Q)}^2 \\[1mm]
     && \hspace*{1cm}  \leq \, - \langle Q_{h_t}^0 P_{h_x}^1 u-u_h,
               \partial_t \overline{v}_h \rangle_{L^2(Q)} +
               \langle \nabla_x (Q_{h_t}^0 P_{h_x}^1 u-u_h),
               \nabla_x \overline{v}_h \rangle_{L^2(Q)} \\[1mm]
     &&  \hspace*{1cm} = \, b_T(Q_{h_t}^0 P_{h_x}^1 u-u_h,\overline{v}_h) \\[1mm]
     && \hspace*{1cm}
        = \, - \langle P_{h_x}^1 u-u,\partial_t \overline{v}_h \rangle_{L^2(Q)}
           + \langle \nabla_x (Q_{h_t}^0 u-u),
           \nabla_x (\overline{v}_h -Q_{h_t}^0\overline{v}_h) \rangle_{L^2(Q)}. 
   \end{eqnarray*}
   Now, using the Cauchy--Schwarz inequality and \eqref{Error Px L2}, the
   first term is estimated by 
   \begin{eqnarray*}
     - \langle P_{h_x}^1 u-u,\partial_t \overline{v}_h \rangle_{L^2(Q)} 
     & \leq &
     \| P_{h_x}^1 u-u \|_{L^2(Q)} \| \partial_t\overline{v}_h \|_{L^2(Q)} \\[1mm]
     & \leq &
     c \, h_x^2 \, |u|_{L^2(0,T;H^2(\Omega))}
     \| \partial_t \overline{v}_h \|_{L^2(Q)}, 
   \end{eqnarray*}
   while for the second term, using integration by parts and
   \eqref{Error Qht}, this gives 
   \begin{eqnarray*}
     \langle \nabla_x (Q_{h_t}^0 u-u),
     \nabla_x (\overline{v}_h -Q_{h_t}^0\overline{v}_h) \rangle_{L^2(Q)} 
     & = & \langle -\Delta_x(Q_{h_t}^0 u-u),
           \overline{v}_h-Q_{h_t}^0 \overline{v}_h \rangle_{L^2(Q)} \\[1mm]
     && \hspace*{-6cm} \leq \, \| \Delta_x(Q_{h_t}^0 u-u) \|_{L^2(Q)}
        \| \overline{v}_h-Q_{h_t}^0 \overline{v}_h \|_{L^2(Q)} \, \leq \,
        c \, h_t^2 \, \| \Delta_x \partial_t u \|_{L^2(Q)}
        \| \partial_t \overline{v}_h \|_{L^2(Q)}. 
   \end{eqnarray*}
   With $h_x\simeq h_t$ and
   $\partial_t \overline{v}_h = Q_{h_t}^0 P_{h_x}^1 u-u_h$, we finally obtain
   \[
     \| Q_{h_t}^0 P_{h_x}^1 u-u_h \|_{L^2(Q)} 
     \, \leq \,
     c \, h_x^2 \, \Big(
     \|u\|_{L^2(0,T;H^2(\Omega))} + \|u\|_{H^1(0,T;H^2(\Omega))} \Big),
   \]
   and the assertion follows.}

\noindent
It remains to consider the case of less regular initial data:

\begin{corollary}
  Let $u \in Y$ and $u_h \in Y_h$ be the unique solutions of the variational
  formulations \eqref{Heat adjoint VF} and \eqref{VF adjoint heat},
  respectively. Assume $ f \equiv 0$, and
  $u_0 \in \widetilde{H}^s(\Omega) := [L^2(\Omega),H^1_0(\Omega)]_{|s}$ for
  some $s \in [0,1]$. Then,
  $u \in H^{1+s,(1+s)/2}(Q) \cap H^{s/2}(0,T;H^1_0(\Omega))$. The error
  estimate \eqref{Error adjoint L2} then gives
  \begin{eqnarray*}
    \| u - u_h \|_{L^2(Q)}
    & \leq & c \, \Big[ h_x^{1+s} \,
             \| u \|_{L^2(0,T;H^{1+s}(\Omega))} +
             h_t^{(1+s)/2} \, \| u \|_{H^{(1+s)/2}(0,T;L^2(\Omega))} \\
    && \hspace*{4cm} + \,
    h_x^{-1} \, h_t \, h_t^{s/2} \, \| u \|_{H^{s/2}(0,T;H^1_0(\Omega))} \Big],
  \end{eqnarray*}
  and in particular for $h_t \simeq h_x^2$ we obtain
  \begin{eqnarray*}
    \| u - u_h \|_{L^2(Q)} & & \\
    && \hspace*{-2cm} \leq \, 
    c \, h_x^{1+s} \, \Big[ \| u \|_{L^2(0,T;H^{1+s}(\Omega))} +
    \| u \|_{H^{(1+s)/2}(0,T;L^2(\Omega))} +
    \| u \|_{H^{s/2}(0,T;H^1_0(\Omega))} \Big] .
  \end{eqnarray*}
  Note that for a discontinuous but piecewise smooth initial datum we have
  $u_0 \in \widetilde{H}^s(\Omega)$, $ s < \frac{1}{2}$.
\end{corollary}

\noindent
It remains to take a closer look on the discrete variational formulation
\eqref{VF adjoint heat}. When writing
\[
  u_h(x,t) \, = \, \sum\limits_{i=1}^{N_t} \sum\limits_{k=1}^{M_x} u_k^i \,
  \psi_i^0(t) \, \varphi_k(x) \, = \,
  \sum\limits_{k=1}^{M_x} u_k(t) \, \varphi_k(x),
\]
and using as test function $v_h(x,t) = \psi_j^1(t) \varphi_\ell(x)$,
$j=1,\ldots,N_t-1$, $\ell=1,\ldots,M_x$, \eqref{VF adjoint heat}
is equivalent to, recall $\psi_j^1(0)=0$ for $j=1,\ldots,N_t-1$,
\begin{eqnarray*}
  && - \sum\limits_{k=1}^{M_x}
  \int_{t_{j-1}}^{t_{j+1}} 
  u_k(t) \, \partial_t \psi_j^1(t) \, dt
  \int_\Omega \varphi_k(x) \varphi_\ell(x) dx \\
  && \hspace*{10mm}
  + \sum\limits_{k=1}^{M_x} \int_{t_{j-1}}^{t_{j+1}}
  u_k(t) \psi_j^1(t) \, dt
  \int_\Omega \nabla_x \varphi_k(x) \cdot \nabla_x \varphi_\ell(x) \, dx \\
  &&  \hspace*{5cm}
     = \, \int_{t_{j-1}}^{t_{j+1}} \int_\Omega f(x,t) \varphi_\ell(x) \, dx \,
  \psi_j^1(t) \, dt, 
\end{eqnarray*}
which we can write as Crank--Nicolson scheme
\begin{equation}\label{adjoint CN}
  M_h (\underline{u}^{j+1}-\underline{u}^j) +
  \frac{1}{2} \, h_t \, K_h (\underline{u}^{j+1} + \underline{u}^j)
  \, = \, \underline{f}^j ,
\end{equation}
with the spatial mass and stiffness matrices $M_h$ and $K_h$. 
Note that the coefficient vector $\underline{u}^j$ describes the solution
which is constant in the time interval $(t_j,t_{j+1})$.
For the remaining test functions
$v_h(x,t) = \psi_0^1(t) \varphi_\ell(x)$, $\ell=1,\ldots,M_x$, we have
\begin{eqnarray*}
  && - \sum\limits_{k=1}^{M_x}
  \int_{t_0}^{t_1} u_k(t) \, \partial_t \psi_0^1(t) \, dt
     \int_\Omega \varphi_k(x) \varphi_\ell(x) dx \\
  && \hspace*{10mm}
  + \sum\limits_{k=1}^{M_x} \int_{t_0}^{t_1} u_k(t) \psi_0^1(t) \,
  dt \int_\Omega \nabla_x \varphi_k(x) \cdot \nabla_x \varphi_\ell(x)
  \, dx \\
  && \hspace*{3cm}
    = \, \int_{t_0}^{t_1} \int_\Omega f(x,t) \varphi_\ell(x) \, dx \,
  \psi_0^1(t) \, dt + \int_\Omega u_0(x) \varphi_\ell(x) \, dx, 
\end{eqnarray*}
i.e., we obtain the implicit Euler scheme

\begin{equation}\label{adjoint implicit Euler}
  M_h \underline{u}^1  + \frac{1}{2} \, h_t \, K_h \underline{u}^1
  = \underline{f}^0 + \underline{u}^0 .
\end{equation}
In particular, the adjoint formulation corresponds to a time stepping scheme,
where for the first step implicit Euler is used with step size $h_t/2$, and
then Crank--Nicolson is applied with step size $h_t$. Such a scheme was already
proposed by Rannacher \cite{Rannacher:1985} known as \emph{Rannacher smoothing},
to regularize the initial datum $u_0$ and get a more stable numerical
approximation. Indeed, one could also see the first implicit Euler step as a
smoothing projection $\Pi_{h_x}^1: L^2(\Omega)\to V_{h_x}$ given as 
\begin{equation}\label{eq:definition smoothing ritz projection}
  \langle \Pi_{h_x}^1 u_0 , v_{h_x} \rangle_{L^2(\Omega)} +
  \frac{h_t}{2} \, \langle \nabla_x \Pi_{h_x}^1 u_0 , \nabla_x v_{h_x}
  \rangle_{L^2(\Omega)} \, = \,
  \langle u_0 , v_{h_x} \rangle_{L^2(\Omega)},\quad \mbox{for all} \;
  v_{h_x}\in V_{h_x}.  
\end{equation}
Then the adjoint formulation corresponds to a perturbed primal formulation in
which the initial datum $u_0$ is projected by $\Pi_{h_x}^1$ rather than by
$Q_{h_x}^1$. Since the projection $\Pi_{h_x}^1$ regularizes the initial datum,
the error analysis of the primal formulation for smooth data can be applied.
To obtain the convergence rate of this scheme, it remains to estimate the
projection error of the initial condition.

\begin{lemma}
  For $u_0\in L^2(\Omega)$ there holds 
  \begin{equation}\label{eq: projection stability for u in L2}
    \| \Pi_{h_x}^1 u_0 \|_{L^2(\Omega)} \, \leq \, \| u_0 \|_{L^2(\Omega)},
    \quad \text{and} \quad 
    \| \nabla_x \Pi_{h_x}^1 u_0 \|_{L^2(\Omega)} \, \leq \,
    \sqrt{2} \, h_t^{-1/2} \, \| u_0 \|_{L^2(\Omega)}, 
  \end{equation}
  as well as 
  \begin{equation}\label{eq: i- projection stability for u in L2}
    \| u_0 - \Pi_{h_x}^1 u_0 \|_{L^2(\Omega)} \, \leq \, \| u_0 \|_{L^2(\Omega)}. 
  \end{equation}
  Moreover, if $u_0\in H^1_0(\Omega)$, then 
  \begin{equation}\label{eq:projection stabiity for u in H1}
    \| \nabla_x \Pi_{h_x}^1 u_0 \|_{L^2(\Omega)} \, \leq \,
    c \, \Big( 1 + h_x^2 \, h_t^{-1} \Big)^{1/2} \,
    \| \nabla_x u_0\|_{L^2(\Omega)},
  \end{equation}
  \begin{equation}\label{eq:projection convergence for u in H1}
    \| u_0-\Pi_{h_x}^1 u_0 \|_{L^2(\Omega)} \, \leq \,
    c \, \Big( h_t + h_x^2 \Big)^{1/2} \, \| \nabla_x u_0 \|_{L^2(\Omega)},
  \end{equation}
  and 
  \begin{equation}\label{eq:projection convergence for u in H^1 in H^1 norm}
    \| \nabla_x(u_0-\Pi_{h_x}^1u_0) \|_{L^2(\Omega)} \, \leq \,
    c \, \Big( 1 + h_x^2 \, h_t^{-1} \Big)^{1/2} \,
    \| \nabla_x u_0 \|_{L^2(\Omega)}.
  \end{equation}
  Moreover, if $u_0\in H^2(\Omega) \cap H^1_0(\Omega)$ is satisfied, then 
  \begin{equation}\label{eq:projection converge for u in H^2}
    \| u_0-\Pi_{h_x}^1 u_0 \|_{L^2(\Omega)} \, \leq \, c \,
    \Big( h_t + h_x^2 \Big) \, | u_0 |_{H^2(\Omega)},
  \end{equation}
  and 
  \begin{equation}\label{eq:projection convergence for u in H^2 H^1 norm}
    \| \nabla_x ( u_0-\Pi_{h_x}^1 u_0) \|_{L^2(\Omega)} \, \leq \,
    c \, \Big( h_t + h_x^4 \, h_t^{-1} + h_x^2 \Big)^{1/2} \,
    | u_0 |_{H^2(\Omega)}. 
  \end{equation}
\end{lemma}

\myproof{When testing \eqref{eq:definition smoothing ritz projection} with
  $v_{h_x} = \Pi_{h_x}^1 u_0 \in V_{h_x}$ this gives
  \[
    \| \Pi_{h_x}^1 u_0 \|_{L^2(\Omega)}^2
    +
    \frac{h_t}{2} \, \| \nabla_x \Pi_{h_x}^1 u_0 \|_{L^2(\Omega)}^2 \, = \,
    \langle u_0 , \Pi_{h_x}^1 u_0 \rangle_{L^2(\Omega)} \, \leq \, 
    \| u_0 \|_{L^2(\Omega)} \| \Pi_{h_x}^1 u_0 \|_{L^2(\Omega)},
  \]
  from which \eqref{eq: projection stability for u in L2} follows.
  Moreover, we compute, using \eqref{eq:definition smoothing ritz projection}
  for $v_h = \Pi_{h_x}^1 u_0$,
  \begin{eqnarray*}
    && \| u_0 - \Pi_{h_x}^1 u_0 \|_{L^2(\Omega)}^2 +
    \frac{h_t}{2} \, \| \nabla_x \Pi_{h_x}^1 u_0 \|_{L^2(\Omega)}^2 \\
    && \hspace*{1cm}
       = \, \| u_0 - \Pi_{h_x}^1 u_0 \|_{L^2(\Omega)}^2 + \frac{h_t}{2} \,
       \langle \nabla_x \Pi_{h_x}^1 u_0 , \nabla_x \Pi_{h_x}^1 u_0
       \rangle_{L^2(\Omega)} \\
    && \hspace*{1cm} = \, \| u_0 - \Pi_{h_x}^1 u_0 \|_{L^2(\Omega)}^2 + 
       \langle u_0 - \Pi_{h_x}^1 u_0 , \Pi_{h_x}^1 u_0 \rangle_{L^2(\Omega)} \\
    && \hspace*{1cm}
       = \, \langle u_0 - \Pi_{h_x}^1 u_0 , u_0 \rangle_{L^2(\Omega)} \, \leq \,
       \| u_0 - \Pi_{h_x}^1 u_0 \|_{L^2(\Omega)} \| u_0 \|_{L^2(\Omega)},
  \end{eqnarray*}
  from which \eqref{eq: i- projection stability for u in L2} follows.
  For $u_0 \in H^1_0(\Omega)$ we can consider the  $H^1_0$ projection
  $P_{h_x}^1 u_0 \in V_{h_x}$, and using the same computations as above as well
  as H\"older's inequality we conclude
   \begin{eqnarray*}
     \| u_0 - \Pi_{h_x}^1 u_0 \|_{L^2(\Omega)}^2 +
     \frac{h_t}{2} \, \| \nabla_x \Pi_{h_x}^1 u_0 \|_{L^2(\Omega)}^2
     & = & \langle u_0 - \Pi_{h_x}^1 u_0 , u_0 \rangle_{L^2(\Omega)} \\
     && \hspace*{-6.5cm} = \,
        \langle u_0 - \Pi_{h_x}^1 u_0 , u_0 - P_{h_x}^1 u_0 \rangle_{L^2(\Omega)}
        +
        \langle u_0 - \Pi_{h_x}^1 u_0 , P_{h_x}^1 u_0 \rangle_{L^2(\Omega)} \\
     && \hspace*{-6.5cm} = \,
        \langle u_0 - \Pi_{h_x}^1 u_0 , u_0 - P_{h_x}^1 u_0 \rangle_{L^2(\Omega)}
        + \frac{h_t}{2} \,
        \langle \nabla_x \Pi_{h_x}^1 u_0 ,
        \nabla_x P_{h_x}^1 u_0 \rangle_{L^2(\Omega)} \\
     && \hspace*{-6.5cm} \leq \,
        \| u_0 - \Pi_{h_x}^1 u_0 \|_{L^2(\Omega)}
        \| u_0 - P_{h_x}^1 u_0 \|_{L^2(\Omega)} +
        \frac{h_t}{2} \, \| \nabla_x \Pi_{h_x}^1 u_0 \|_{L^2(\Omega)}
        \| \nabla_x P_{h_x}^1 u_0 \|_{L^2(\Omega)} \\
     && \hspace*{-6.5cm} \leq \, \left(
        \| u_0 - \Pi_{h_x}^1 u_0 \|_{L^2(\Omega)}^2 +
        \frac{h_t}{2} \, \| \nabla_x \Pi_{h_x}^1 u_0 \|^2_{L^2(\Omega)}
        \right)^{1/2} \\
     && \hspace*{-3cm} \cdot \left(
        \| u_0 - P_{h_x}^1 u_0 \|_{L^2(\Omega)}^2 +
        \frac{h_t}{2} \, \| \nabla_x P_{h_x}^1 u_0 \|^2_{L^2(\Omega)}
        \right)^{1/2},
  \end{eqnarray*}
  i.e.,
  \[
    \| u_0 - \Pi_{h_x}^1 u_0 \|_{L^2(\Omega)}^2 +
    \frac{h_t}{2} \, \| \nabla_x \Pi_{h_x}^1 u_0 \|^2_{L^2(\Omega)}
    \, \leq \,
    \| u_0 - P_{h_x}^1 u_0 \|_{L^2(\Omega)}^2 +
    \frac{h_t}{2} \, \| \nabla_x P_{h_x}^1 u_0 \|^2_{L^2(\Omega)} .
  \]
  With the error estimate \eqref{Error Px L2} and the stability of 
  $P_{h_x}^1 : H^1_0(\Omega) \to V_{h_x} \subset H^1_0(\Omega)$
  we further obtain
  \[
    \| u_0 - P_{h_x}^1 u_0 \|_{L^2(\Omega)}^2
    +
    \frac{h_t}{2} \, \| \nabla_x \Pi_{h_x}^1 u_0 \|^2_{L^2(\Omega)}
    \, \leq \,
    c \, \Big[ h_x^2 + h_t \Big] \, \| \nabla_x u_0 \|_{L^2(\Omega)}^2,
  \]
  which gives \eqref{eq:projection stabiity for u in H1} and
  \eqref{eq:projection convergence for u in H1}. Moreover,
  \eqref{eq:definition smoothing ritz projection} is equivalent to
  \[
    \frac{h_t}{2} \, \langle \nabla_x (\Pi_{h_x}^1 u_0 - P_{h_x}^1u_0),
    \nabla_x v_h \rangle_{L^2(\Omega)}
    =
    \langle u_0 - \Pi_{h_x}^1 u_0 , v_h \rangle_{L^2(\Omega)} -
    \frac{h_t}{2} \, \langle \nabla_x P_{h_x}^1u_0,
    \nabla_x v_h \rangle_{L^2(\Omega)}
  \]
  for all $v_h \in V_{h_x}$. In particular for
  $v_h = \Pi_{h_x}^1u_0 - P_{h_x}^1u_0$ this gives
  \begin{eqnarray}
    \frac{h_t}{2} \, \nonumber
    \| \nabla_x (\Pi_{h_x}^1 u_0 - P_{h_x}^1 u_0) \|_{L^2(\Omega)}^2
    && \\
    && \hspace*{-40mm}
       = \, \langle u_0 - \Pi_{h_x}^1 u_0 , \Pi_{h_x}^1 u_0 - P_{h_x}^1 u_0
       \rangle_{L^2(\Omega)} \label{Proof Pihx Step 1} \\
    && \hspace*{-20mm} + \nonumber
       \frac{h_t}{2} \, \langle \nabla_x P_{h_x}^1u_0,
       \nabla_x (P_{h_x}^1 u_0 - \Pi_{h_x}^1 u_0) \rangle_{L^2(\Omega)} \\
    && \hspace*{-40mm} = \, \| u_0 - \Pi_{h_x}^1 u_0 \|_{L^2(\Omega)} \nonumber
       \| \Pi_{h_x}^1 u_0 - P_{h_x}^1 u_0 \|_{L^2(\Omega)} \\
    && \hspace*{-20mm} + \nonumber
       \frac{h_t}{2} \, \| \nabla_x P_{h_x}^1u_0 \|_{L^2(\Omega)}
       \| \nabla_x (P_{h_x}^1 u_0 - \Pi_{h_x}^1 u_0) \|_{L^2(\Omega)} \\
    && \hspace*{-40mm} \leq \, \| u_0 - \Pi_{h_x}^1 u_0 \|_{L^2(\Omega)}
       \Big[ \| \Pi_{h_x}^1 u_0 - u_0 \|_{L^2(\Omega)} + \nonumber
       \| u_0 - P_{h_x}^1 u_0 \|_{L^2(\Omega)} \Big] \\
    && \hspace*{-20mm} + \, \nonumber
       \frac{h_t}{2} \, \| \nabla_x P_{h_x}^1u_0 \|_{L^2(\Omega)}
       \| \nabla_x (P_{h_x}^1 u_0 - \Pi_{h_x}^1 u_0) \|_{L^2(\Omega)} \\
    && \hspace*{-40mm} \leq \, c \, \Big( h_t + h_x^2 \Big)^{1/2} \, \nonumber
       \| \nabla_x u_0 \|_{L^2(\Omega)} \\
    && \nonumber \hspace*{-10mm} \cdot \Big[
       c \, \Big( h_t + h_x^2 \Big)^{1/2} \,
       \| \nabla_x u_0 \|_{L^2(\Omega)}
       + c \, h_x \, \| \nabla_x u_0 \|_{L^2(\Omega)} \Big] \\
    && \hspace*{-20mm} + \, \nonumber
       \frac{h_t}{2} \, \| \nabla_x u_0 \|_{L^2(\Omega)}
       \| \nabla_x (P_{h_x}^1 u_0 - \Pi_{h_x}^1 u_0) \|_{L^2(\Omega)} \\
    && \hspace*{-40mm} \leq \, c \, \Big( h_t + h_x^2 \Big) \, \nonumber
       \| \nabla_x u_0 \|_{L^2(\Omega)}^2 + 
       \frac{h_t}{2} \, \| \nabla_x u_0 \|_{L^2(\Omega)}
       \| \nabla_x (P_{h_x}^1 u_0 - \Pi_{h_x}^1 u_0) \|_{L^2(\Omega)},
  \end{eqnarray}
  when using the error estimates \eqref{eq:projection convergence for u in H1}
  and \eqref{Error Px L2}, and the stability estimate \eqref{Stability Ph}.
  Hence we have shown
  \begin{eqnarray*}
    && \| \nabla_x (\Pi_{h_x}^1 u_0 - P_{h_x}^1 u_0) \|_{L^2(\Omega)}^2 \\
    && \hspace*{2mm}
       \leq \, 2 \, c \, \Big( 1 + h_x^2 \, h_t^{-1} \Big) \,
       \| \nabla_x u_0 \|_{L^2(\Omega)}^2 + 
       \| \nabla_x u_0 \|_{L^2(\Omega)}
       \| \nabla_x (P_{h_x}^1 u_0 - \Pi_{h_x}^1 u_0) \|_{L^2(\Omega)},
  \end{eqnarray*}
  from which
  \[
    \| \nabla_x (\Pi_{h_x}^1 u_0 - P_{h_x}^1 u_0) \|_{L^2(\Omega)} \, \leq \,
    c \, \Big( 1 + h_x^2 \, h_t^{-1} \Big)^{1/2} \,
    \| \nabla_x u_0 \|_{L^2(\Omega)}
  \]
  follows. With
  \[
    \| \nabla_x (u_0 - \Pi_{h_x}^1u_0) \|_{L^2(\Omega)} \, \leq \,
    \| \nabla_x (u_0 - P_{h_x}^1u_0) \|_{L^2(\Omega)} +
    \| \nabla_x (P_{h_x}^1 u_0 - \Pi_{h_x}^1u_0) \|_{L^2(\Omega)}
  \]
  and using \eqref{Error Px H1} for $s=1$
  we conclude \eqref{eq:projection convergence for u in H^1 in H^1 norm}.

  It remains to consider the case $u_0 \in H^2(\Omega) \cap H^1_0(\Omega)$.
  Note that \eqref{Proof Pihx Step 1} can be rewritten as
  \begin{eqnarray*}
    && \| u_0 - \Pi_{h_x}^1 u_0 \|^2_{L^2(\Omega)}
       + 
       \frac{h_t}{2} \, 
       \| \nabla_x (\Pi_{h_x}^1 u_0 - P_{h_x}^1 u_0) \|_{L^2(\Omega)}^2 \\
    && \hspace*{2mm} = \,
       \langle u_0 - \Pi_{h_x}^1 u_0 , u_0 - P_{h_x}^1 u_0 \rangle_{L^2(\Omega)}
       + \frac{h_t}{2} \, \langle \nabla_x P_{h_x}^1u_0,
       \nabla_x (P_{h_x}^1 u_0 - \Pi_{h_x}^1 u_0) \rangle_{L^2(\Omega)} \\
    && \hspace*{2mm} = \,
       \langle u_0 - \Pi_{h_x}^1 u_0 , u_0 - P_{h_x}^1 u_0 \rangle_{L^2(\Omega)}
       + \frac{h_t}{2} \, \langle \nabla_x u_0,
       \nabla_x (P_{h_x}^1 u_0 - \Pi_{h_x}^1 u_0) \rangle_{L^2(\Omega)},
  \end{eqnarray*}
  using the definition \eqref{Def Phx} of $P_{h_x}^1$. Using integration
  by parts, we further conclude
  \begin{eqnarray*}
    && \| u_0 - \Pi_{h_x}^1 u_0 \|^2_{L^2(\Omega)}
       + 
       \frac{h_t}{2} \, 
       \| \nabla_x (\Pi_{h_x}^1 u_0 - P_{h_x}^1 u_0) \|_{L^2(\Omega)}^2 \\
    && \hspace*{1cm} = \,
       \langle u_0 - \Pi_{h_x}^1 u_0 , u_0 - P_{h_x}^1 u_0 \rangle_{L^2(\Omega)}
       + \frac{h_t}{2} \, \langle - \Delta_x u_0,
       P_{h_x}^1 u_0 - \Pi_{h_x}^1 u_0 \rangle_{L^2(\Omega)} \\
    && \hspace*{1cm} \leq \,
       \| u_0 - \Pi_{h_x}^1 u_0 \|_{L^2(\Omega)}
       \| u_0 - P_{h_x}^1 u_0 \|_{L^2(\Omega)} \\
    && \hspace*{2.5cm}
       + \frac{h_t}{2} \, \| \Delta_x u_0 \|_{L^2(\Omega)}
       \| P_{h_x}^1 u_0 - \Pi_{h_x}^1 u_0 \|_{L^2(\Omega)} \\
    && \hspace*{1cm} \leq \,
       \| u_0 - \Pi_{h_x}^1 u_0 \|_{L^2(\Omega)}
       \| u_0 - P_{h_x}^1 u_0 \|_{L^2(\Omega)} \\
    && \hspace*{2.5cm}
       + \frac{h_t}{2} \, \| \Delta_x u_0 \|_{L^2(\Omega)} \Big[
       \| P_{h_x}^1 u_0 - u_0 \|_{L^2(\Omega)} +
       \| u_0 - \Pi_{h_x}^1 u_0 \|_{L^2(\Omega)} \Big] \\
    && \hspace*{1cm} \leq \, c \, h_x^2 \,
       \| u_0 - \Pi_{h_x}^1 u_0 \|_{L^2(\Omega)}
       | u_0 |_{H^2(\Omega)} \\
    && \hspace*{2.5cm}
       + \frac{h_t}{2} \, | u_0 |_{H^2(\Omega)} \Big[ c \, h_x^2 \,
       | u_0 |_{H^2(\Omega)} +
       \| u_0 - \Pi_{h_x}^1 u_0 \|_{L^2(\Omega)} \Big],
  \end{eqnarray*}
  where we have used the error estimate \eqref{Error Px L2} for $s=2$.
  Hence we have shown
  \[
    \| u_0 - \Pi_{h_x}^1 u_0 \|^2_{L^2(\Omega)}
    \, \leq \,
    c \, \Big( h_x^2 + h_t \Big)
    | u_0 |_{H^2(\Omega)} \| u_0 - \Pi_{h_x}^1 u_0 \|_{L^2(\Omega)} 
    + c \, h_x^2 \, h_t \, | u_0 |_{H^2(\Omega)}^2,
  \]
  from which
  \[
    \| u_0 - \Pi_{h_x}^1 u_0 \|_{L^2(\Omega)}
    \leq
    c \, \Big( h_x^2 + h_t \Big) |u_0|_{H^2(\Omega)}
  \]
  follows, i.e., \eqref{eq:projection converge for u in H^2}.
  With this we further conclude
  \begin{eqnarray*}
    && \| \nabla_x (\Pi_{h_x}^1 u_0 - P_{h_x}^1 u_0) \|_{L^2(\Omega)}^2 \\
    && \hspace*{1cm} \leq \, c \, \Big( 1 + h_x^2 \, h_t^{-1} \Big) 
             \| u_0 - \Pi_{h_x}^1 u_0 \|_{L^2(\Omega)} | u_0 |_{H^2(\Omega)} 
             + c \, h_x^2 \, | u_0 |_{H^2(\Omega)}^2 \\
    && \hspace*{1cm} \leq \, c \, \Big[
             \Big( 1 + h_x^2 \, h_t^{-1} \Big) \Big( h_x^2 + h_t \Big)
             + h_x^2 \Big]
                | u_0 |_{H^2(\Omega)}^2 .
  \end{eqnarray*}
  With the triangle inequality
  \begin{eqnarray*}
    && \| \nabla_x (u_0 - \Pi_{h_x}^1 u_0) \|_{L^2(\Omega)}^2 \\
    && \hspace*{1cm} \leq \, 2 \, \| \nabla_x (u_0 - P_{h_x}^1 u_0) \|_{L^2(\Omega)}^2
             + 2 \, \| \nabla_x (P_{h_x}^1 u_0 - \Pi_{h_x}^1 u_0) \|_{L^2(\Omega)}^2 \\
    && \hspace*{1cm} \leq \, c \, \Big[
             \Big( 1 + h_x^2 \, h_t^{-1} \Big) \Big( h_x^2 + h_t \Big)
             + h_x^2 \Big] | u_0 |_{H^2(\Omega)}^2
  \end{eqnarray*}
  we finally conclude \eqref{eq:projection convergence for u in H^2 H^1 norm}.}

\noindent
When considering the error estimates 
\eqref{eq: i- projection stability for u in L2} and
\eqref{eq:projection convergence for u in H1} we can use a space interpolation
argument to conclude the following result:

\begin{theorem}\label{thm:error estimate ritz projection}
  Assume $u_0 \in \widetilde{H}^s(\Omega) := [L^2(\Omega),H^1_0(\Omega)]_{|s}$
  for some $s \in [0,1]$. Then there holds the error estimate
  \begin{equation}\label{eq:error estimate ritz projection}
    \| u_0 - \Pi_{h_x}^1 u_0 \|_{L^2(\Omega)} \, \leq \,
    c \, \Big( h_x^2 + h_t \Big)^{s/2} \,
    \| u_0 \|_{\widetilde{H}^s(\Omega)} .
  \end{equation}
\end{theorem}

\section{Numerical results}\label{sec:numerics}
As a numerical example we consider the heat equation \eqref{Heat equation}
in the one-dimensional spatial domain $\Omega = (0,1)$ and $T=1$.
The space-time finite element discretization of the primal variational
formulation \eqref{Heat equation FEM} is done by using the ansatz space
$X_{0,h} = W_{h_t}^1 \otimes V_{h_x}$ and the test space
$Y_h = W_{h_t}^0 \otimes V_{h_x}$, while for the adjoint formulation
\eqref{VF adjoint heat} we use $Y_h$ as ansatz space, and
$X_{T,h} = \overline{W}_{h_t}^1 \otimes V_{h_x}$ as test space. First,
we consider the smooth solution
\begin{equation}\label{smooth solution heat}
u(x,t) = \cos (\pi t) \sin (\pi x) 
\end{equation}
where the given data $f$ and $u_0$ are computed from $u$. The numerical
results for both the primal and the adjoint formulation using $h_t=h_x$
are given in Table \ref{Table heat smooth}, where we observe first order
convergence in $Y=L^2(0,T;H^1_0(\Omega))$ in both cases, as predicted by
Theorem \ref{Theorem error primal smooth X} and
Theorem \ref{thm:Error adjoint Y smooth}, respectively. However, we only
observe second order convergence in $L^2(Q)$ for the primal formulation
as explained by Lemma \ref{Lemma L2 Error primal smooth}, while the error for the adjoint approach converges linearly, see
Corollary \ref{cor:Error L2 adjoint smooth}. To obtain second
order convergence in $L^2(Q)$ also for the adjoint formulation, we need to
postprocess the solution, see Theorem \ref{thm:error reconstruction in time},
or use parabolic scaling, see Table~\ref{Table heat smooth parabolic}, and
Table~\ref{Table heat smooth uniform reconstructed}.

\begin{table}[ht]
  \begin{center}
\begin{tabular}{rrlclc}
  \hline
  $N_x$ & $N_t$
  & $\| u - u_h \|_{L^2(Q)}$ & eoc
  & $\| u - u_h \|_Y$ & eoc  \\
  \hline
  \multicolumn{6}{c}{primal formulation} \\ 
  \hline
    4 &   4 & $4.137 \cdot 10^{-2}$ &      & $3.597 \cdot 10^{-1}$ &      \\
    8 &   8 & $1.027 \cdot 10^{-2}$ & 2.01 & $1.784 \cdot 10^{-1}$ & 1.01 \\
   16 &  16 & $2.562 \cdot 10^{-3}$ & 2.00 & $8.907 \cdot 10^{-2}$ & 1.00 \\
   32 &  32 & $6.403 \cdot 10^{-4}$ & 2.00 & $4.452 \cdot 10^{-2}$ & 1.00 \\
   64 &  64 & $1.601 \cdot 10^{-4}$ & 2.00 & $2.226 \cdot 10^{-2}$ & 1.00 \\
  128 & 128 & $4.001 \cdot 10^{-5}$ & 2.00 & $1.113 \cdot 10^{-2}$ & 1.00 \\
  256 & 256 & $1.000 \cdot 10^{-5}$ & 2.00 & $5.565 \cdot 10^{-3}$ & 1.00 \\
  \hline
  \multicolumn{6}{c}{adjoint formulation} \\
  \hline
    4 &   4 & $1.129 \cdot 10^{-1}$ &      & $4.988 \cdot 10^{-1}$ &      \\
    8 &   8 & $5.661 \cdot 10^{-2}$ & 1.00 & $2.512 \cdot 10^{-1}$ & 0.99 \\
   16 &  16 & $2.833 \cdot 10^{-2}$ & 1.00 & $1.258 \cdot 10^{-1}$ & 1.00 \\
   32 &  32 & $1.417 \cdot 10^{-2}$ & 1.00 & $6.295 \cdot 10^{-2}$ & 1.00 \\
   64 &  64 & $7.085 \cdot 10^{-3}$ & 1.00 & $3.148 \cdot 10^{-2}$ & 1.00 \\
  128 & 128 & $3.543 \cdot 10^{-3}$ & 1.00 & $1.574 \cdot 10^{-2}$ & 1.00 \\
  256 & 256 & $1.771 \cdot 10^{-3}$ & 1.00 & $7.870 \cdot 10^{-3}$ & 1.00 \\
  \hline
\end{tabular}
\end{center}
\caption{Convergence results for the primal and adjoint formulations in the
  case of the smooth solution \eqref{smooth solution heat}, $h_t = h_x$.}
\label{Table heat smooth}
\end{table}
\begin{table}[ht]
\begin{center}
\begin{tabular}{rrlclc}
  \hline
  \multicolumn{6}{c}{adjoint formulation} \\
  \hline
  $N_x$ & $N_t$
  & $\| u -u_h \|_{L^2(Q)}$ & eoc
  & $\| u - u_h \|_Y$ & eoc \\
  \hline
    4 &    16 & $5.110 \cdot 10^{-2}$ &      & $4.784 \cdot 10^{-1}$ &      \\
    8 &    64 & $1.054 \cdot 10^{-2}$ & 2.28 & $2.044 \cdot 10^{-1}$ & 1.23 \\
   16 &   256 & $2.452 \cdot 10^{-3}$ & 2.10 & $9.509 \cdot 10^{-2}$ & 1.10 \\
   32 &  1024 & $5.944 \cdot 10^{-4}$ & 2.04 & $4.597 \cdot 10^{-2}$ & 1.05 \\
   64 &  4096 & $1.465 \cdot 10^{-4}$ & 2.02 & $2.261 \cdot 10^{-2}$ & 1.02 \\
  128 & 16384 & $3.638 \cdot 10^{-5}$ & 2.01 & $1.122 \cdot 10^{-2}$ & 1.01 \\
  \hline
\end{tabular}
\end{center}
\caption{Convergence results for the adjoint formulation in the case of the
  smooth solution \eqref{smooth solution heat} with parabolic scaling,
  $h_t = h_x^2$.}
\label{Table heat smooth parabolic}
\end{table}
\begin{table}[ht]
\begin{center}
  \begin{tabular}{rrlclc}
  \hline
  \multicolumn{6}{c}{adjoint formulation} \\
  \hline
  $N_x$ & $N_t$ & $\| u - \mathcal{R}_{h_t} u_h \|_{L^2(Q)}$ & eoc
  & $\| u - \mathcal{R}_{h_t}u_h \|_Y$ & eoc \\
   \hline
    4 &   4 & $6.186 \cdot 10^{-2}$ &       & $3.847 \cdot 10^{-1}$ &       \\
    8 &   8 & $1.779 \cdot 10^{-2}$ & 1.798 & $1.822 \cdot 10^{-1}$ & 1.078 \\
   16 &  16 & $4.705 \cdot 10^{-3}$ & 1.918 & $8.957 \cdot 10^{-2}$ & 1.024 \\
   32 &  32 & $1.205 \cdot 10^{-3}$ & 1.965 & $4.459 \cdot 10^{-2}$ & 1.006 \\
   64 &  64 & $3.047 \cdot 10^{-4}$ & 1.984 & $2.227 \cdot 10^{-2}$ & 1.002 \\
  128 & 128 & $7.656 \cdot 10^{-5}$ & 1.992 & $1.113 \cdot 10^{-2}$ & 1.000 \\
  256 & 256 & $1.919 \cdot 10^{-5}$ & 1.996 & $5.565 \cdot 10^{-3}$ & 1.000 \\
  \hline
\end{tabular}
\end{center}
\caption{Convergence results for the adjoint formulation with reconstruction
    in the case of the smooth solution \eqref{smooth solution heat}
    and uniform scaling, $h_t=h_x$.}
\label{Table heat smooth uniform reconstructed}
\end{table}

\noindent
As a second example, we consider, similar to \cite{Osterby2003}, the solution
\begin{equation}\label{singular solution}
  u(x,t)
  =
  \sum_{\ell=1}^{\infty}
  \frac{4}{\ell\pi}
  \sin\left(\frac{\ell\pi}{2}\right)
  \sin\left(\frac{\ell\pi}{4}\right)
  e^{-\ell^2\pi^2 t}
  \sin(\ell\pi x).
\end{equation}
which is a solution of the homogeneous heat equation \eqref{Heat equation}
with $f \equiv 0$, and with the initial datum
\[
  u_0(x) = \left\{
    \begin{array}{ll}
      1, & \mbox{if} \; x\in(0.25,0.75), \\[2mm]
      0, & \mbox{else}.
    \end{array} 
  \right.
\]
The numerical results for both the
primal and adjoint formulations are given in Table \ref{Table heat singular}. 
In the case of orthotropic scaling, $h_t = h_x$, we observe divergence 
of the primal formulation in $Y$, as can be explained by
\eqref{Cea primal full X}, where $c_{S,h}\simeq h_x$, and using the best
approximation in $X$. For the adjoint formulation, though, we see a rate of
$1/4$, as predicted by Theorem \ref{thm:error estimate ritz projection} for all $s<1/2$.
The rates for the error in $L^2(Q)$ are not yet explained by the theory presented. 
When using parabolic scaling, $h_t=h_x^2$, in Table \ref{Table heat singular parabolic} in both cases we obtain
comparable results with orders of convergence which are due to the
reduced regularity of the solution as explained by \eqref{Cea primal full X}
and \eqref{adjoint Cea} for the primal and adjoint formulation, respectively,
using $c_{S,h}\simeq 1$, and the best approximation in the corresponding norms.

\begin{figure}
  \centering
  \begin{tikzpicture}
    \begin{axis}[
      width=10cm,
      height=6cm,
      xlabel={$x$},
      ylabel={$u(x,\,t=0.1)$},
      legend pos=outer north east,
      grid=major,
    ]
      \addplot[
        thick,
        color=blue,
        mark=o,
        mark options={solid},
      ]
        table[x=x, y=uprim, col sep=comma] {img/solution_values32.csv};
      \addlegendentry{primal}
  
      \addplot[
        thick,
        color=red,
        mark=square*,
        mark options={solid},
      ]
        table[x=x, y=udual, col sep=comma] {img/solution_values32.csv};
      \addlegendentry{adjoint}
  
      \addplot[
        thick,
        color=black,
        dashed,
      ]
        table[x=x, y=uex, col sep=comma] {img/solution_values32.csv};
      \addlegendentry{analytical}
    \end{axis}
  \end{tikzpicture}
  \caption{Comparison of the primal and adjoint approach at $t=0.1$.}
  \label{fig:discontinuity}
\end{figure}
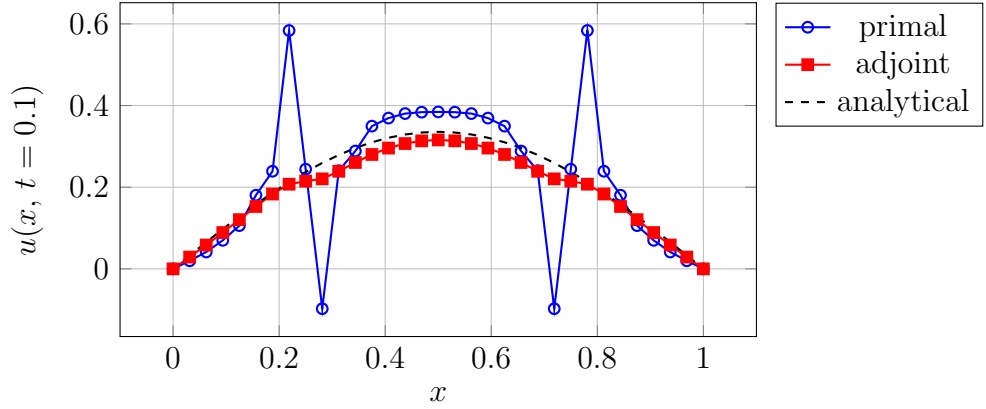

\begin{table}[ht]
  \begin{center}
\begin{tabular}{rrlclc}
  \hline
  $N_x$ & $N_t$
  & $\|u-u_h\|_{L^2(Q)}$ & eoc
  & $\|u-u_h\|_Y$ & eoc \\
  \hline
  \multicolumn{6}{c}{primal formulation} \\
  \hline
   4 &  4 & $8.878 \cdot 10^{-2}$ &      & $8.377 \cdot 10^{-1}$ &       \\
   8 &  8 & $7.222 \cdot 10^{-2}$ & 0.30 & $1.115 \cdot 10^{0}$  & -0.41 \\
  16 & 16 & $5.545 \cdot 10^{-2}$ & 0.38 & $1.599 \cdot 10^{0}$  & -0.52 \\
  32 & 32 & $4.091 \cdot 10^{-2}$ & 0.44 & $2.269 \cdot 10^{0}$  & -0.50 \\
  64 & 64 & $2.959 \cdot 10^{-2}$ & 0.47 & $3.183 \cdot 10^{0}$  & -0.49 \\
  \hline
  \multicolumn{6}{c}{adjoint formulation} \\
  \hline
   4 &  4 & $8.944 \cdot 10^{-2}$ &      & $3.498 \cdot 10^{-1}$ &      \\
   8 &  8 & $5.273 \cdot 10^{-2}$ & 0.76 & $2.753 \cdot 10^{-1}$ & 0.35 \\
  16 & 16 & $3.073 \cdot 10^{-2}$ & 0.78 & $2.356 \cdot 10^{-1}$ & 0.22 \\
  32 & 32 & $1.832 \cdot 10^{-2}$ & 0.75 & $2.025 \cdot 10^{-1}$ & 0.22 \\
  64 & 64 & $1.091 \cdot 10^{-2}$ & 0.75 & $1.725 \cdot 10^{-1}$ & 0.23 \\
  \hline
\end{tabular}
\end{center}
\caption{Convergence results for the primal and adjoint formulations for the
  singular solution \eqref{singular solution} with discontinuous
  initial data, $h_t=h_x$.}
\label{Table heat singular}
\end{table}
\begin{table}[ht]
  \begin{center}
\begin{tabular}{rrlclc}
  \hline
  $N_x$ & $N_t$
  & $\|u-u_h\|_{L^2(Q)}$ & eoc
  & $\|u-u_h\|_Y$ & eoc \\
  \hline
  \multicolumn{6}{c}{primal formulation} \\
  \hline
   4 &   16 & $2.265 \cdot 10^{-2}$ &      & $3.091 \cdot 10^{-1}$ & \\
   8 &   64 & $7.467 \cdot 10^{-3}$ & 1.60 & $1.979 \cdot 10^{-1}$ & 0.64 \\
  16 &  256 & $2.642 \cdot 10^{-3}$ & 1.50 & $1.397 \cdot 10^{-1}$ & 0.50 \\
  32 & 1024 & $9.314 \cdot 10^{-4}$ & 1.50 & $9.854 \cdot 10^{-2}$ & 0.50 \\
  64 & 4096 & $3.264 \cdot 10^{-4}$ & 1.51 & $6.920 \cdot 10^{-2}$ & 0.51 \\
  \hline
  \multicolumn{6}{c}{adjoint formulation} \\
  \hline
   4 &   16 & $3.073 \cdot 10^{-2}$ &      & $2.435 \cdot 10^{-1}$ &      \\
   8 &   64 & $1.088 \cdot 10^{-2}$ & 1.50 & $1.726 \cdot 10^{-1}$ & 0.50 \\
  16 &  256 & $3.848 \cdot 10^{-3}$ & 1.50 & $1.221 \cdot 10^{-1}$ & 0.50 \\
  32 & 1024 & $1.360 \cdot 10^{-3}$ & 1.50 & $8.629 \cdot 10^{-2}$ & 0.50 \\
  64 & 4096 & $4.807 \cdot 10^{-4}$ & 1.50 & $6.090 \cdot 10^{-2}$ & 0.50 \\
  \hline
\end{tabular}
\end{center}
\caption{Convergence results for the primal and adjoint formulations for the
  singular solution \eqref{singular solution} with discontinuous
  initial data and parabolic scaling, $h_t=h_x^2$.}
\label{Table heat singular parabolic}
\end{table}
\noindent
Even though the convergence behavior is similar, the adjoint approach does
not overshoot at the discontinuities of the initial data. This behavior is
exemplified in Figure \ref{fig:discontinuity} for the primal and adjoint
approach for $N_x=N_t=32$, see also \cite{Osterby2003} for a more detailed
discussion of the Crank--Nicolson oscillations.

\noindent
Although we have presented almost all error estimates in the case of lowest
order finite element spaces only, these results can be generalized to higher
order polynomial ansatz spaces. To ensure optimal convergence, a related higher
regularity of the solution is required, excluding discontinuous initial data.
For the smooth solution as given in \eqref{smooth solution heat} we provide
some numerical results, in which we first consider second order finite elements in
space, and second order finite elements in time as ansatz functions, but
piecewise linear ones for testing. In the adjoint formulation we switch ansatz
and test spaces accordingly. The numerical results as given in
Table \ref{Table heat smooth p2} indicate optimal orders of convergence, as
expected. Finally, in Table \ref{Table heat smooth p3} we provide related
results for cubic finite elements in space and time, and second order test
functions in time for the primal formulation. For the adjoint
formulation, second order functions in time are used as ansatz, while
third order finite elements are used for testing.
Again we observe optimal orders of convergence, as expected.

\begin{table}[ht]
  \begin{center}
\begin{tabular}{rrlclc}
  \hline
  $N_x$ & $N_t$ & $\| u - u_h \|_{L^2(Q)}$ & eoc & $\| u - u_h \|_Y$ & eoc \\
  \hline
  \multicolumn{6}{c}{primal formulation} \\
  \hline
  2 &  2 & $5.78 \cdot 10^{-2}$ &      & $1.12 \cdot 10^{-1}$ & \\
  4 &  4 & $1.71 \cdot 10^{-2}$ & 1.76 & $2.02 \cdot 10^{-1}$ & -0.85 \\
  8 &  8 & $2.49 \cdot 10^{-3}$ & 2.78 & $5.13 \cdot 10^{-2}$ &  1.98 \\
 16 & 16 & $3.25 \cdot 10^{-4}$ & 2.94 & $1.28 \cdot 10^{-2}$ &  2.00 \\
 32 & 32 & $4.10 \cdot 10^{-5}$ & 2.98 & $3.20 \cdot 10^{-3}$ &  2.00 \\
 64 & 64 & $5.15 \cdot 10^{-6}$ & 2.99 & $7.99 \cdot 10^{-4}$ &  2.00 \\
  \hline
  \multicolumn{6}{c}{adjoint formulation} \\
  \hline
   2 &  2 & $4.64 \cdot 10^{-2}$ &      & $2.05 \cdot 10^{-1}$ & \\
   4 &  4 & $1.20 \cdot 10^{-2}$ & 1.95 & $5.17 \cdot 10^{-2}$ & 1.98 \\
   8 &  8 & $2.92 \cdot 10^{-3}$ & 2.04 & $1.28 \cdot 10^{-2}$ & 2.01 \\
  16 & 16 & $7.22 \cdot 10^{-4}$ & 2.02 & $3.20 \cdot 10^{-3}$ & 2.00 \\
  32 & 32 & $1.80 \cdot 10^{-4}$ & 2.01 & $7.98 \cdot 10^{-4}$ & 2.00 \\
  64 & 64 & $4.49 \cdot 10^{-5}$ & 2.00 & $2.00 \cdot 10^{-4}$ & 2.00 \\
  \hline
\end{tabular}
\end{center}
\caption{Convergence results for the primal and adjoint formulations using
  quadratic/linear finite elements in time, and quadratic/quadratic finite elements
  in space in the case of the smooth solution \eqref{smooth solution heat}, $h_t=h_x$.}
\label{Table heat smooth p2}
\end{table}

\begin{table}[ht]
  \begin{center}
\begin{tabular}{rrlclc}
  \hline
  $N_x$ & $N_t$ & $\| u - u_h \|_{L^2(Q)}$ & eoc & $\| u - u_h \|_Y$ & eoc \\
  \hline
  \multicolumn{6}{c}{primal formulation} \\
  \hline
  2 &  2 & $4.48 \cdot 10^{-2}$ &      & $7.76 \cdot 10^{-2}$ & \\
  4 &  4 & $2.20 \cdot 10^{-3}$ & 4.35 & $6.86 \cdot 10^{-3}$ & 3.50 \\
  8 &  8 & $1.31 \cdot 10^{-4}$ & 4.07 & $7.25 \cdot 10^{-4}$ & 3.24 \\
 16 & 16 & $8.01 \cdot 10^{-6}$ & 4.03 & $8.39 \cdot 10^{-5}$ & 3.11 \\
 32 & 32 & $4.98 \cdot 10^{-7}$ & 4.01 & $1.02 \cdot 10^{-5}$ & 3.03 \\
 64 & 64 & $3.11 \cdot 10^{-8}$ & 4.00 & $1.27 \cdot 10^{-6}$ & 3.01 \\
  \hline
  \multicolumn{6}{c}{adjoint formulation} \\
  \hline
   2 &  2 & $6.47 \cdot 10^{-3}$ &      & $1.99 \cdot 10^{-2}$ & \\
   4 &  4 & $7.78 \cdot 10^{-4}$ & 3.06 & $2.43 \cdot 10^{-3}$ & 3.03 \\
   8 &  8 & $9.59 \cdot 10^{-5}$ & 3.02 & $3.01 \cdot 10^{-4}$ & 3.01 \\
  16 & 16 & $1.19 \cdot 10^{-5}$ & 3.01 & $3.75 \cdot 10^{-5}$ & 3.00 \\
  32 & 32 & $1.49 \cdot 10^{-6}$ & 3.00 & $4.68 \cdot 10^{-6}$ & 3.00 \\
  64 & 64 & $1.86 \cdot 10^{-7}$ & 3.00 & $5.85 \cdot 10^{-7}$ & 3.00 \\
  \hline
\end{tabular}
\end{center}
\caption{Convergence results for the primal and adjoint formulations using
  cubic/quadratic finite elements in time and cubic/cubic finite elements in
  space in the case of the smooth solution \eqref{smooth solution heat}, $h_t=h_x$.}
\label{Table heat smooth p3}
\end{table}

\section{Conclusions}\label{sec:Conclusions}
In this paper, we analyzed primal and adjoint space-time variational
formulations for parabolic evolution equations using tensor product
discretizations in space and time. For lowest order ansatz and test
spaces, the primal formulation recovers the Crank--Nicolson method.
At the same time, its space-time interpretation reveals a conditional
stability property, where the stability constant depends on a CFL-type
relation between the temporal and spatial mesh sizes.

This conditional stability becomes relevant in the error analysis for
solutions of low regularity, in particular for discontinuous or nonsmooth
initial data. The numerical experiments confirm the deteriorated behaviour
of the primal method in this regime. To incorporate the initial condition
differently, we derived an adjoint space-time formulation when applying
integration by parts in time. For regular data, this formulation yields
comparable convergence rates after a suitable postprocessing in time. For
rough data, however, it shows an improved behaviour and naturally leads
to a Rannacher-type smoothing of the initial data.

The Galerkin--Petrov space-time formulation therefore provides an
operator-theoretic framework, based on inf-sup stability, for analyzing
both classical time stepping schemes and their behaviour for nonsmooth data.
For low-regularity problems, the adjoint formulation is preferable, whereas
for smooth solutions the primal formulation remains competitive. Although
the detailed analysis in this work was carried out for lowest order
approximations, the underlying arguments extend to higher order test and
ansatz spaces, as can be observed in the numerical examples. 

\bigskip

\noindent
\textbf{Acknowledgement:}
This work is supported by the joint
DFG/FWF Collaborative Research Centre CREATOR
(DFG: Project-ID 492661287/TRR361; FWF: 10.55776/F90) at
TU Darmstadt, TU Graz and JKU Linz.

\begin{appendix}

    \section{Regularity estimates}
    In this appendix we summarize some well known results on the
    regularity of solutions for the initial boundary value problem
    \eqref{Heat equation}, see, e.g., \cite[Sect. 7.1.3]{Evans:1998}.
    Since the heat equation is linear, we first
    consider the homogeneous heat equation, and afterwards the inhomogeneous
    one but with zero initial conditions. We assume that $\Omega \subset {\mathbb{R}}^n$ is bounded with smooth
    boundary $\partial \Omega$, or polygonal/polyhedral but convex.

    \subsection{Homogeneous heat equation}
    We first consider the homogeneous heat equation with homogeneous
    Dirichlet boundary conditions, but some initial conditions,
    \begin{equation}\label{Appendix:homogeneous heat}
      \partial_t u - \Delta_x u = 0 \; \mbox{in} \; Q = \Omega \times (0,T),
      \quad
      u = 0 \; \mbox{on} \; \Sigma = \partial \Omega \times (0,T),
      \quad
      u(0)=u_0 \; \mbox{in} \; \Omega .
    \end{equation}
    The solution of (\ref{Appendix:homogeneous heat}) is given by the
    Fourier series
    \begin{equation}\label{Appendix:homogeneous heat solution}
      u(x,t) = \sum\limits_{k=1}^\infty u_k \,
      e^{- \lambda_k t} \phi_k(x), \quad
      u_k = \int_\Omega u_0(x) \, \phi_k(x) \, dx ,
    \end{equation}
    where $\phi_k \in H^1_0(\Omega)$ are the normalized Dirichlet
    eigenfunctions of the Laplacian,
    \[
      - \Delta \phi_k = \lambda_k \phi_k \quad \mbox{in} \; \Omega, \quad
      \phi_k = 0 \quad \mbox{on} \; \partial \Omega, \quad
      \int_\Omega \phi_k(x) \, \phi_j(x) \, dx = \delta_{kj} \, .
    \]
    Based on \eqref{Appendix:homogeneous heat solution} and depending
    on the regularity of the given initial datum $u_0$ we can state the
    following results:

    \begin{lemma}\label{Appendix:homogeneous heat regularity}
      Let $u$ as given in \eqref{Appendix:homogeneous heat solution}
      be the solution of the initial boundary value problem
      \eqref{Appendix:homogeneous heat}. Depending on the regularity
      of the initial datum $u_0$ there holds:
      \begin{enumerate}
      \item For $u_0 \in L^2(\Omega)$ we have
        $u \in L^2(0,T;H^1_0(\Omega)) \cap H^1(0,T;H^{-1}(\Omega))$,
        satisfying
        \begin{equation}\label{Appendix:homogeneous heat Regularity L2 H11}
          \| u \|^2_{L^2(0,T;H^1_0(\Omega))} + \| u \|^2_{H^1(0,T;H^{-1}(\Omega))}
          \, \leq \, \| u_0 \|^2_{L^2(\Omega)} ,
        \end{equation}
        as well as $u \in H^{1,1/2}(Q) :=
        L^2(0,T;H^1_0(\Omega)) \cap H^{1/2}(0,T;L^2(\Omega))$,
        satisfying
        \begin{equation}\label{Appendix:homogeneous heat Regularity L2 H112}
          \| u \|^2_{L^2(0,T;H^1_0(\Omega))} + \| u \|^2_{H^{1/2}(0,T;L^2(\Omega))}
          \, \leq \, c \, \| u_0 \|^2_{L^2(\Omega)} .
        \end{equation}
      \item For $u_0 \in H^1_0(\Omega)$,
        $u \in H^{2,1}(Q) := L^2(0,T;H^2(\Omega)) \cap H^1(0,T;L^2(\Omega))$,
        satisfying
        \begin{equation}\label{Appendix:homogeneous heat Regularity L2 H21}
          \| u \|^2_{L^2(0,T;H^2(\Omega))} + \| u \|^2_{H^1(0,T;L^2(\Omega))}
          \, \leq \, \| \nabla_x u_0 \|^2_{L^2(\Omega)} .
        \end{equation}
      \item For $u_0 \in H^1_0(\Omega) \cap H^2(\Omega)$, there holds
        \begin{equation}\label{Appendix:homogeneous heat Regularity uxt}
          \| \nabla_x \partial_t u \|^2_{L^2(Q)} \, \leq \, \frac{1}{2} \,
          |u_0|^2_{H^2(\Omega)},
        \end{equation}
        and
        \begin{equation}\label{Appendix:homogeneous heat Regularity utt}
          \| u \|_{H^2(0,T;H^{-1}(\Omega))}^2 \, \leq \, \frac{1}{2} \,
          \| u_0 \|_{H^2(\Omega)}^2 .
        \end{equation}       
      \end{enumerate}
    \end{lemma}

    \noindent
    When using space interpolation arguments, we derive regularity
    results in an intermediate scale of Sobolev spaces. 
    
    \begin{corollary}
      Assume $u_0 \in \widetilde{H}^s(\Omega) :=
      [L^2(\Omega),H^1_0(\Omega)]_{|s}$ for some $s \in [0,1]$. From
      \eqref{Appendix:homogeneous heat Regularity L2 H112} and
      \eqref{Appendix:homogeneous heat Regularity L2 H21} 
      we then conclude $u \in H^{1+s,(1+s)/2}(Q)$, and
      \begin{equation}\label{Appendix:homogeneous heat u0 Hss2}
        \| u \|_{L^2(0,T;H^{1+s}(\Omega))}^2 +
        \| u \|_{H^{(1+s)/2}(0,T;L^2(\Omega))}^2 \, \leq \, c \,
        \| u_0 \|^2_{\widetilde{H}^s(\Omega)} .
      \end{equation}
      Moreover, using \eqref{Appendix:homogeneous heat Regularity L2 H112}
      and \eqref{Appendix:homogeneous heat Regularity uxt} we obtain
      $u \in H^{s/2}(0,T;H^1_0(\Omega))$, satisfying
      \begin{equation}\label{Appendix:homogeneous heat Regularity Hs2}
        \| u \|_{H^{s/2}(0,T;H^1_0(\Omega))} \, \leq \, c \,
        \| u_0 \|_{\widetilde{H}^s(\Omega)} .
      \end{equation}
    \end{corollary}

    \noindent
    Finally, from \eqref{Appendix:homogeneous heat Regularity L2 H11}
    and \eqref{Appendix:homogeneous heat Regularity utt} we conclude
    the following estimate:

    \begin{corollary}
      For the solution $u$ of \eqref{Appendix:homogeneous heat} with
      $u_0 \in H^1_0(\Omega)$ there holds
      \begin{equation}\label{Appendix:homogeneous heat u32}
        \| u \|^2_{H^{3/2}(0,T;H^{-1}(\Omega))} \, \leq \, c \,
        \| u_0 \|^2_{H^1(\Omega)} .
      \end{equation}
    \end{corollary}

    \subsection{Inhomogeneous heat equation}
    Next we consider the inhomogeneous heat equation with zero Dirichlet
    boundary conditions, and zero initial conditions,
    \begin{equation}\label{Appendix:Heat f}
      \partial_t u - \Delta_x u = f \; \mbox{in} \; Q, \quad
      u = 0 \; \mbox{on} \; \Sigma, \quad
      u(0)=0 \; \mbox{in} \; \Omega .
    \end{equation}
    The solution of \eqref{Appendix:Heat f} is given by the Fourier series
    \begin{equation}\label{Appendix:Heat f solution}
      u(x,t) = \sum\limits_{k=1}^\infty u_k(t) \varphi_k(x),
    \end{equation}
    where the coefficients $u_k(t)$ are solutions of the ordinary
    differential equations
    \[
      u_k'(t) + \lambda_k u_k(t) \, = \, f_k(t) \, := \,
      \langle f(t),\varphi_k \rangle_\Omega \quad \mbox{for} \; t \in (0,T),
      \quad u_k(0)=0 .
    \]

    \begin{lemma}
      Let $u$ as given in \eqref{Appendix:Heat f solution} be the solution of
      the initial boundary value problem \eqref{Appendix:Heat f}.
      Depending on the regularity of $f$ there holds:
      \begin{enumerate}
      \item For $f \in L^2(0,T;H^{-1}(\Omega))$, we have
        $u \in L^2(0,T;H^1_0(\Omega)) \cap H^1_{0,}(0,T;H^{-1}(\Omega))$,
        with
        \begin{equation}\label{Appendix:Heat f Regularity H11}
          \| u \|^2_{L^2(0,T;H^1_0(\Omega))} \, + \,
          \| u \|_{H^1(0,T;H^{-1}(\Omega))} ^2
          \, \leq \, 2 \, \| f \|_{L^2(0,T;H^{-1}(\Omega))}^2 ,
        \end{equation}
        as well as
        $ u \in H^{1,1/2}_{0;0,}(Q) :=
        L^2(0,T;H^1_0(\Omega)) \cap H_{0,}^{1/2}(0,T;L^2(\Omega))$,
        satisfying
        \begin{equation}\label{Appendix:Heat f Regularity H112}
          \| u \|^2_{L^2(0,T;H^1_0(\Omega))} +
          \| u \|^2_{H^{1/2}_{0,}(0,T;L^2(\Omega))} \, \leq \,
          2 \, \| f \|_{L^2(0,T;H^{-1}(\Omega))}^2 .
        \end{equation}
      \item For $f \in L^2(Q)$,
        $u \in L^2(0,T;H^1_0(\Omega) \cap H^2(\Omega)) \cap
        H^1_{0,}(0,T;L^2(\Omega))$, satisfying
        \begin{equation}\label{Appendix:Heat f Regularity H21}
          \| u \|^2_{L^2(0,T;H^1_0(\Omega) \cap H^2(\Omega))} +
          \| u \|^2_{H^1_{0,}(0,T;L^2(\Omega))} \, \leq \, 2 \,
          \| f \|^2_{L^2(Q)} .  
        \end{equation}
      \item For $f \in L^2(0,T;H^1_0(\Omega))$, there holds
        \begin{equation}\label{Appendix:Heat f Regularity uxt}
          \| \nabla_x \partial_t u \|_{L^2(Q)}^2 \, \leq \,
          \| f \|_{L^2(0,T;H^1_0(\Omega))}^2 .  
        \end{equation}
      \item If $f \in H^1_{0,}(0,T;H^{-1}(\Omega))$ is satisfied, then
        \begin{equation}\label{Appendix:Heat f Regularity utt}
          \| \partial_{tt} u \|_{L^2(0,T;H^{-1}(\Omega))} \, \leq \,
          \sqrt{2} \, \| f \|_{H^1(0,T;H^{-1}(\Omega))} .
        \end{equation}
      \end{enumerate}
    \end{lemma}

    \noindent
    Again we can use space interpolation arguments to formulate
    regularity results in an intermediate scale of Sobolev spaces.
    From \eqref{Appendix:Heat f Regularity H112} and
    \eqref{Appendix:Heat f Regularity H21} we have:
    
    \begin{corollary}
      Assume $f \in L^2(0,T;H^{-1+s}(\Omega)) :=
      [L^2(0,T;H^{1-s}_0(\Omega))]^*$ for some $s \in [0,1]$.
      Then, we conclude $u \in H^{1+s,(1+s)/2}(Q)$, and
      \begin{equation}\label{Appendix:Heat f Regularity Hss/2}
        \| u \|^2_{L^2(0,T;H^{1+s}(\Omega))} +
        \| u \|^2_{H^{(1+s)/2}(0,T;L^2(\Omega))} \, \leq \, c \,
        \| f \|^2_{L^2(0,T;H^{-1+s}(\Omega))} \, .
      \end{equation}  
    \end{corollary}
    
    \noindent
    Finally, from \eqref{Appendix:Heat f Regularity H11} and
    \eqref{Appendix:Heat f Regularity utt} we conclude the following
    estimate:

    \begin{corollary}
      For the solution $u$ of \eqref{Appendix:Heat f} with
      $f \in H^{1/2}_{0,}(0,T;H^{-1}(\Omega))$ there holds
      \begin{equation}\label{Appendix:Heat f Regularity u32}
        \| u \|^2_{H^{3/2}(0,T;H^{-1}(\Omega))} \, \leq \,
        c \, \| f \|^2_{H^{1/2}(0,T;H^{-1}(\Omega))} .
      \end{equation}
    \end{corollary}

    \section{Equivalence of inf-sup stability constants}
    Because for any bounded, linear isomorphism $B : X \to Y^*$ it holds
    that $(B^{-1})^* = (B^*)^{-1}$, and
    $\| B \|_{X \to Y^*} = \| B^* \|_{Y \to X^*}$, the continuous inf-sup
    constants for $B$ and $B^*$ are the same, i.e., 
    \[
      \inf_{0 \neq u \in X} \sup_{0\neq v\in Y}
      \frac{\langle Bu,v \rangle_{Y^*\times Y}}{\|u\|_X \|v\|_Y} 
      \, = \,
      \inf_{0 \neq v\in Y} \sup_{0\neq u\in X}
      \frac{\langle B^* v,u \rangle_{X^* \times X}}{\| u \|_X \| v \|_Y}. 
    \]
    The same conclusion holds true for the discrete inf-sup constant,
    e.g., \cite[Proposition 3.4.3]{BoffiBrezziFortin:2013}. For completeness,
    and because our discrete trial and test spaces differ slightly, we
    include a proof of this result in the next lemma.

    \begin{lemma}\label{lem:discrete inf-sup equality}
      It holds that 
      \[
        \inf_{0 \neq u_h\in X_{0,h}} \sup_{0\neq v_h\in Y_h}
        \frac{b(u_h,v_h)}{\|u_h\|_X\|v_h\|_Y} \, = \,
        \inf_{0 \neq v_h \in Y_h} \sup_{0\neq u_h\in X_{T,h}}
        \frac{b_T(v_h,u_h)}{\| v_h \|_Y \| u_h \|_X},   
      \]
      where $b(\cdot,\cdot) : X_{0,h}\times Y_h \to \mathbb{R}$ and
      $b_T(\cdot,\cdot):Y_h\times X_{T,h}\to\mathbb{R}$ are defined
      in \eqref{eq:bilinearform b} and \eqref{Heat adjoint VF}, respectively. 
    \end{lemma}

    \myproof{We first note that the time reversal map
      $\iota_T v(t) := v(T-t)$, when solely acting on time, is an
      isometric isomorphism as a map $\iota_T:Y_h\to Y_h$, and as a
      map $\iota_T:X_{0,h}\to X_{T,h}$. Moreover,
      $\partial_t \iota_T = -\iota_T \partial_t$, $\iota_T^2 = \text{id}$,
      and the $L^2$ adjoint satisfies $\iota_T^* = \iota_T$. Hence,
      for all $u \in X_0$ and $v \in Y$ we have that
      \begin{eqnarray*}
        b(u,v)
        & = & \int_0^T \int_\Omega \Big[
              \partial_t u(x,t) \, v(x,t) +
              \nabla_x u(x,t) \cdot \nabla_x v(x,t) \Big] \, dt \, dx \\
        & = & \int_0^T \int_\Omega \Big[
              \partial_t u(x,t) \, \iota_T^2v(x,t) +
              \nabla_x u(x,t) \cdot \nabla_x \iota_T^2 v(x,t) \Big] dt \, dx \\
        & = & \int_0^T \int_\Omega \Big[ \iota_T^* \partial_t u(x,t) \,
              \iota_T v(x,t) + \nabla_x \iota_T^* u(x,t) \cdot \nabla_x
              \iota_T v(x,t) \Big] dt \, dx \\
        & = & \int_0^T \int_\Omega \Big[
              \iota_T \partial_t u(x,t) \, \iota_T v(x,t) +
              \nabla_x \iota_T u(x,t) \cdot \nabla_x \iota_T v(x,t)
              \Big] dt \, dx \\
        & = & \int_0^T \int_\Omega \Big[ - \partial_t \iota_T u(x,t) \,
              \iota_T v(x,t) + \nabla_x \iota_T u(x,t) \cdot
              \nabla_x \iota_T v(x,t) \Big] dt \, dx \\
        & = & b_T(\iota_T v,\iota_T u).
      \end{eqnarray*}
      Thus, for the operators $B_{0,h} : X_{0,h} \to Y_h^*$, and
      $B_{T,h} : Y_h \to X_{T,h}^*$, defined as 
      \[
        \langle B_{0,h}u_h,v_h \rangle_{Y^* \times Y}
        \, = \,
        b(u_h,v_h), \quad
        \langle B_{T,h} v_h,u_h \rangle_{X^* \times X} \, = \, b_T(v_h,u_h),
      \]
      there holds the relation 
      \begin{eqnarray*}
        \langle B_{0,h} u_h,v_h \rangle_{Y^* \times Y}
        & = & b(u_h,v_h) \, = \, b_T(\iota_T v_h,\iota_T u_h) \\
        & = & \langle B_{T,h} \iota_T v_h, \iota_T u_h \rangle_{X^* \times X}
              \, = \, \langle \iota_T^* B_{T,h}^* \iota_T u_h,v_h
              \rangle_{Y^* \times Y}
      \end{eqnarray*}
      for all $u_h\in X_{0,h}$, and $v_h \in Y_h$, i.e.,
      $B_{0,h} = \iota_T^* B_{T,h}^* \iota_T$. Now, we compute the discrete
      inf-sup constant for the primal problem as
      \begin{eqnarray*}
        \inf_{0 \neq u_h\in X_{0,h}} \sup_{0\neq v_h \in Y_h}
        \frac{b(u_h,v_h)}{\|u_h\|_X \|v_h\|_Y}
        & = &
        \inf_{0 \neq  u_h\in X_{0,h}} \sup_{0\neq v_h \in Y_h}
        \frac{\langle B_{0,h} u_h,v_h \rangle_{Y^* \times Y}}
              {\| u_h \|_X \| v_h \|_Y} \\
        & = &
        \inf_{0 \neq u_h \in X_{0,h}}
        \frac{\| B_{0,h} u_h \|_{Y^*}}{\| u_h \|_X}.   
      \end{eqnarray*}
      Moreover,  
      \begin{eqnarray*}
        \inf_{0 \neq u_h \in X_{0,h}}
        \frac{\| B_{0,h} u_h \|_{Y^*}}{\|u_h\|_X}
        & = & \inf_{0 \neq v_h \in Y_h^*}
              \frac{\| v_h \|_{Y^*}}{\| B_{0,h}^{-1}v_h \|_X} \\
        & = & \left( \sup_{0\neq v_h\in Y_h^*}
              \frac{\| B_{0,h}^{-1} v_h\|_X}{\| v_h\|_{Y^*}}\right)^{-1}
              \, = \, \| B_{0,h}^{-1} \|^{-1}_{Y_h^* \to X_{0,h}}.
      \end{eqnarray*}
      Using that $\| B_{0,h}^{-1} \|_{Y_h^* \to X_{0,h}} =
      \| B_{0,h}^{-1,*} \|_{X_{0,h}^* \to Y_h}$, and that
      $B_{0,h}^{-1,*} = B_{0,h}^{\ast,-1} =
      (\iota_T^* B_{T,h} \iota_T)^{-1}$, we then compute
      \begin{eqnarray*}
        \| B_{0,h}^{-1} \|^{-1}_{Y_h^* \to X_{0,h}}
        & = & \| (\iota_T^* B_{T,h} \iota_T)^{-1} \|_{X_{0,h}^* \to Y_h}^{-1} \\
        & = & \inf_{0 \neq v_h\in Y_h} \sup_{0\neq u_h\in X_{0,h}}
              \frac{\langle B_{T,h} \iota_T v_h,\iota_T u_h
              \rangle_{X^* \times X}}{\| u_h \|_X \| v_h \|_Y} \\
        & = & \inf_{0 \neq v_h\in Y_h} \sup_{0\neq u_h\in X_{0,h}}
              \frac{b_T(\iota_T v_h,\iota_T u_h)}
              {\|u_h\|_X \|v_h\|_Y} \\
        & = & \inf_{0 \neq v_h\in Y_h} \sup_{0\neq u_h\in X_{T,h}}
              \frac{b_T(v_h,u_h)}{\| u_h \|_X \| v_h \|_Y}. 
      \end{eqnarray*}
      This concludes the proof.}
    
\end{appendix}

\end{document}